\documentclass[preprint,12pt]{elsarticle}

\usepackage[hidelinks]{hyperref}
\usepackage{booktabs}
\usepackage{tabularx}
\usepackage{seqsplit}
\usepackage{amsmath,amssymb,amsfonts}
\usepackage{bm}
\usepackage{geometry}
\usepackage{microtype}
\usepackage{soul}
\usepackage{array}
\usepackage{makecell}
\usepackage{multirow}
\usepackage{enumitem}

\usepackage{lineno}

\usepackage{color,xcolor}
\usepackage{soul}

\journal{Powder Technology}

\begin{document}

\begin{frontmatter}
    \title{Geometry-Aware Latent Representation Learning for Drag Coefficient Prediction of Irregular Particles}

    \author[label1]{Zihao Li}
    \author[label2]{Shota Deguchi}
    \author[label3]{Mitsuteru Asai}
    \author[label4]{Shunhua Chen}
    \author[label2]{Naoto Mitsume}
    \affiliation[label1]{
        organization={Degree Programs in Systems and Information Engineering, University of Tsukuba},
        addressline={1-1-1 Tennodai},
        city={Tsukuba},
        postcode={305-8573},
        state={Ibaraki},
        country={Japan}
    }
    \affiliation[label2]{
        organization={Institute of Systems and Information Engineering, University of Tsukuba},
        addressline={1-1-1 Tennodai},
        city={Tsukuba},
        postcode={305-8573},
        state={Ibaraki},
        country={Japan}
    }
    \affiliation[label3]{
        organization={Department of Civil Engineering, Kyushu University},
        addressline={744 Motooka},
        city={Nishi-ku},
        postcode={819-0395},
        state={Fukuoka},
        country={Japan}
    }
    \affiliation[label4]{
        organization={School of Ocean Engineering and Technology, Sun Yat-sen University},
        addressline={Southern Marine Science and Engineering Guangdong Laboratory (Zhuhai)},
        city={Zhuhai},
        postcode={519082},
        state={Guangdong},
        country={People's Republic of China}
    }

    \newpage
    \begin{abstract}
        Predicting the drag coefficient \(C_\mathrm{d}\) of irregular particles requires representing their complex three-dimensional geometry in a form suitable for predictive models.
        Empirical drag correlations and many data-driven models represent particle geometry through a small set of hand-crafted shape descriptors.
        Geometric differences omitted by these descriptors are unavailable to the predictor.
        We propose a variational autoencoder-multilayer perceptron (VAE-MLP) framework that learns a latent representation directly from voxelized particle geometry.
        Geometry reconstruction and drag coefficient prediction jointly guide the learning of this representation.
        We compare VAE-MLP with a descriptor-based MLP (Descriptor-MLP) and established empirical correlations using CFD data for synthetic particle geometries at Reynolds numbers of 10--200.
        When the evaluated geometric variations are represented in the training data, the two learned models achieve comparable accuracy for regular ellipsoids, while VAE-MLP yields lower errors for irregular particles.
        When high-order surface perturbations are excluded from training, VAE-MLP achieves lower overall mean absolute relative error in all five training runs. At the strongest tested perturbation, the mean errors across these runs are 4.97\% for VAE-MLP and 9.88\% for Descriptor-MLP.
        These results support learning geometric representations directly from particle geometry to improve drag prediction for the irregular geometries and unseen surface perturbations examined here.
    \end{abstract}




\end{frontmatter}


\newpage
\section{Introduction}

The hydrodynamic resistance experienced by irregular particles in fluid flows is commonly characterized by the drag coefficient \(C_{\mathrm{d}}\).
The pressure and viscous stresses that determine this resistance depend sensitively on particle geometry and flow conditions, resulting in a nonlinear dependence of drag coefficient on these factors.
Accurate drag prediction is therefore relevant to a wide range of particulate-flow applications, including fluidized systems~\cite{mahajan2018}, microplastic transport in environmental flows~\cite{goral2023}, and sediment transport involving naturally irregular particles~\cite{cheng2025}.

To approximate this dependence, a conventional approach is to use empirical drag correlations, which provide tractable estimates through explicit functions of a prescribed set of variables.
Some correlations are developed for specific particle families, such as cylindrical, disk-like, cubic, and ellipsoidal particles~\cite{gabitto2008,rajitha2006,chen2015,ke2018}.
Others seek broader applicability across different particle shapes by incorporating compact geometric descriptors, ranging from sphericity-based measures~\cite{haider1989drag} to combinations of sphericity and projected-area measures~\cite{ganser1993rational,holzer2008} and descriptors derived from particle dimensions such as flatness and elongation~\cite{bagheri2016drag}.
Their algebraic simplicity and low computational cost make such correlations attractive as drag closures in Eulerian--Lagrangian simulations~\cite{zhong2016demcfddem,li2024numerical}.

Despite their widespread use, most correlations remain applicable only over specific parameter ranges~\cite{michaelides2023}.
These limitations have motivated alternative predictive approaches for irregular-particle drag.
Data-driven machine-learning methods offer one such approach by learning the nonlinear dependence directly from data, without requiring a pre-specified functional form. Data-driven drag models have employed standalone neural-network regression~\cite{zhu2022,yan2019}, alternative machine-learning regressors and ensemble architectures~\cite{hussain2025,wang2026}, as well as hybrid and physics-informed formulations that incorporate prior modelling knowledge~\cite{presareyes2024,tian2025}.
However, increasing the flexibility of the predictive model does not by itself remove limitations associated with the geometric representation.
Particle geometry is still commonly represented using hand-crafted descriptors, from single measures such as sphericity~\cite{yan2019} to broader sets of selected geometric features~\cite{presareyes2024,tian2025}. Hand-crafted descriptors provide a compact and interpretable representation of particle geometry, but necessarily summarize only selected aspects of a three-dimensional shape. For irregular particles, geometrically distinct shapes can exhibit similar descriptor values while producing different drag responses~\cite{michaelides2023}. This ambiguity raises a representation-level question: whether low-dimensional descriptors retain sufficient predictive information as particle geometry becomes more complex.

Recent studies have begun to move beyond conventional scalar descriptors by using richer engineered geometric features~\cite{hu2024} or direct three-dimensional particle representations~\cite{ashwin2022}.
Latent representation learning offers a complementary approach in which compact geometric features are learned directly from high-dimensional particle geometry.
Variational autoencoders (VAEs), for example, learn such representations by reconstructing the input while regularizing the latent space~\cite{kingma2013autoencoding}. Hwang et al.~\cite{hwang2021} applied this approach to particle-force prediction by first learning geometric features and then using them in a force predictor. Related studies extended this representation learning strategy to other particle-force modeling settings~\cite{hwang2023,hwang2024}. These studies demonstrated the feasibility of learned geometric representations, but did not systematically isolate their predictive value relative to conventional hand-crafted descriptors under controlled changes in geometric complexity and training coverage.

The relative predictive value of learned and hand-crafted geometric representations therefore remains unclear as geometric complexity increases and as test geometries depart from the variations represented during training.
We investigate this question along two controlled dimensions: geometric complexity and the training coverage of geometric variations.
To this end, we propose a jointly trained VAE-MLP framework in which geometry reconstruction and drag-prediction supervision shape the latent representation, and compare its predictive performance with that of a descriptor-based MLP under matched evaluation settings.
The main contributions of this work are as follows:
\begin{enumerate}[label=(\arabic*)]
    \item We develop a jointly trained VAE-MLP framework in which geometry reconstruction, drag prediction, and latent regularization jointly shape a compact geometry-aware representation.

    \item We design controlled experiments spanning regular ellipsoids, overall-shape irregularity, and combined overall-shape and fine-scale surface variations, while separately evaluating geometric variations represented during training and fine-scale variations introduced only at evaluation.

    \item We show that the advantage of the learned representation is most pronounced for unseen fine-scale surface variations and persists across five training seeds. Further analyses show that the two models differ in how they use geometric information to predict drag on these unseen geometries.
\end{enumerate}

\section{Methodology}

\subsection{Problem Formulation}

This study considers the drag coefficient \(C_\mathrm{d}\) of a fixed particle in a uniform flow at a prescribed incident angle. Dynamic effects such as translation, flow-induced rotation, and particle deformation are explicitly excluded. Here, \(\mathcal{G}\) denotes the complete three-dimensional spatial configuration of the particle at the prescribed orientation. Accordingly, different particle orientations relative to the flow are encoded within \(\mathcal{G}\) through rotations of the particle geometry, rather than supplied as an additional scalar input. The surrounding fluid is assumed to be single-phase, incompressible, and Newtonian.

For any given particle configuration, the flow condition is parameterized solely by the Reynolds number \(\mathrm{Re}\), while the corresponding drag response is quantified by the drag coefficient:
\begin{equation}
    \mathrm{Re}
    =
    \frac{\rho_\infty U_\infty D_{\mathrm{eq}}}{\eta}
    =
    \frac{U_\infty D_{\mathrm{eq}}}{\nu},
    \qquad
    \nu=\frac{\eta}{\rho_\infty},
    \label{eq:reynolds_number_definition}
\end{equation}
\begin{equation}
    C_{\mathrm{d}}
    =
    \frac{F_{\mathrm{d}}}
    {\tfrac{1}{2}\rho_\infty U_\infty^2 A_{\mathrm{ref}}},
    \label{eq:drag_coefficient_definition}
\end{equation}
where \(U_\infty\) is the free-stream velocity, \(D_{\mathrm{eq}}\) is the geometry-specific volume-equivalent diameter, \(\rho_\infty\) is the reference fluid density, \(\eta\) and \(\nu\) are the dynamic and kinematic viscosities, respectively, and \(A_{\mathrm{ref}}\) is the drag-reference area. The force \(F_{\mathrm{d}}\) is the sum of the pressure and viscous force components in the incoming-flow direction. The specific definitions of the flow and geometric quantities, including \(D_{\mathrm{eq}}\) and \(A_{\mathrm{ref}}\), are detailed in section~\ref{sec:cfd_reference_data}.

With \(\mathcal{G}\) denoting the particle configuration and \(\mathrm{Re}\) serving as the sole flow-condition input, let \(\mathcal{R}(\mathcal{G})\) denote the geometric representation through which the particle configuration can be used by the predictor. In this study, \(\mathcal{R}(\mathcal{G})\) is instantiated either as a learned latent representation or as a set of selected hand-crafted descriptors. The prediction problem can then be written as
\begin{equation}
    \hat{C}_{\mathrm{d}}
    =
    f_{\boldsymbol{\psi}}
    \left(
        \mathcal{R}(\mathcal{G}),
        \mathrm{Re}
    \right),
    \label{eq:prediction_formulation}
\end{equation}
where \(\hat{C}_{\mathrm d}\) is the predicted drag coefficient and \(f_{\boldsymbol{\psi}}(\cdot)\) denotes the predictive mapping parameterized by \(\boldsymbol{\psi}\).
The geometric representation \(\mathcal{R}({\mathcal G})\) distinguishes the two data-driven approaches considered in this study:
\begin{equation}
    \mathcal{R}(\mathcal{G})
    =
    \begin{cases}
        \left(D_{\mathrm{eq}},\,\Phi,\,\Phi_{\mathrm{cross}}\right),
        & \text{Descriptor-MLP}, \\
        \boldsymbol{z}(\mathcal{G}),
        & \text{VAE-MLP}.
    \end{cases}
    \label{eq:geometry_representation}
\end{equation}
Descriptor-MLP therefore represents the particle geometry using a prescribed set of hand-crafted descriptors, whereas VAE-MLP learns a latent representation directly from the voxelized three-dimensional particle geometry. Details of descriptor-MLP are provided in Section~\ref{sec:comparison_baselines}. The following subsections describe the VAE-MLP framework and its latent-representation learning procedure in detail.

\subsection{Geometry-Aware Latent Representation Learning with VAE-MLP}
\label{sec:overall_framework}

\begin{figure}[htbp]
    \centering
    \includegraphics[width=0.9\linewidth]{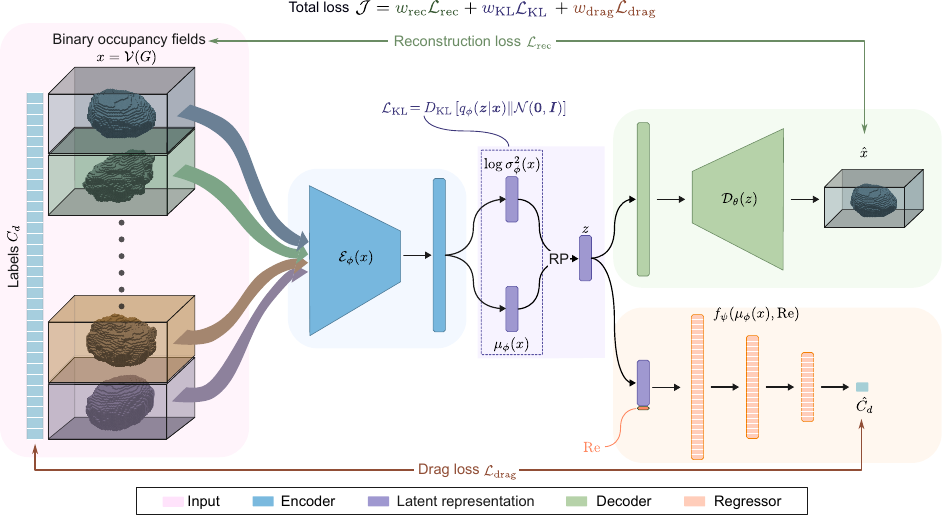}
    \caption{VAE-MLP implementation of geometry-aware latent representation learning for drag coefficient prediction. In the latent representation block, RP indicates the reparameterization trick.}
    \label{fig:methodology_overview}
\end{figure}

The proposed VAE-MLP learns a compact latent representation directly from voxelized particle geometry and uses this learned representation for geometry reconstruction and drag-coefficient prediction.
Figure~\ref{fig:methodology_overview} summarizes the overall architecture and information flow of the proposed VAE-MLP model.
To represent the particle geometry, we first discretize it into a three-dimensional binary voxel occupancy field:
\begin{equation}
    \bm{x}
    =
    \mathcal{V}(\mathcal{G}),
    \qquad
    \bm{x}
    \in
    \{0,1\}^{N_v\times N_v\times N_v}.
    \label{eq:voxelization_mapping}
\end{equation}
The voxelization is adopted as a straightforward way of converting the continuous particle geometry into a regular-grid representation. Details of voxelization are described in Section~\ref{sec:voxelization}.

The variational encoder \(\mathcal E_{\boldsymbol\phi}\) maps the voxelized geometry \(\bm x\) to an approximate posterior distribution in the latent space. During training, a sampled latent vector \(\bm z\) from this distribution is shared by two downstream branches. The decoder \(\mathcal D_{\boldsymbol\theta}\) reconstructs the voxelized geometry, while the drag regressor \(f_{\boldsymbol\psi}\) combines the latent vector \(\bm z\) with \(\mathrm{Re}\) to predict \(C_{\mathrm d}\).
The encoder, decoder, and drag regressor are trained jointly. Consequently, the latent representation is shaped not only to reconstruct particle geometry but also to retain geometric variations relevant to drag prediction, while KL regularization constrains the latent distribution.

Layer-level specifications for the encoder, decoder, and drag regressor, including channel numbers, convolutional parameters, hidden-layer widths, and tensor dimensions, are provided in~\ref{app:architecture_specification}.

\subsection{Variational Latent Encoding}

The variational encoder uses three-dimensional convolutional layers to progressively encode the voxelized particle geometry \(\bm{x}\) across spatial scales. The extracted features are then mapped to two \(d\)-dimensional parameter vectors that define the approximate posterior distribution.

For each input geometry, the encoder produces a posterior mean vector \(\bm{\mu}_{\boldsymbol{\phi}}(\bm{x})\) and a log-variance vector \(\bm{v}_{\boldsymbol{\phi}}(\bm{x})\), where
\begin{equation}
    \bm{v}_{\boldsymbol{\phi}}(\bm{x})
    =
    \log
    \bm{\sigma}_{\boldsymbol{\phi}}^2(\bm{x}).
    \label{eq:posterior_parameters}
\end{equation}
These quantities parameterize a diagonal Gaussian approximate posterior,
\begin{equation}
    q_{\boldsymbol{\phi}}(\bm{z}\mid\bm{x})
    =
    \mathcal{N}\!\left(
        \bm{z};
        \bm{\mu}_{\boldsymbol{\phi}}(\bm{x}),
        \operatorname{diag}\!\left[
            \exp\bigl(\bm{v}_{\boldsymbol{\phi}}(\bm{x})\bigr)
        \right]
    \right).
    \label{eq:variational_posterior}
\end{equation}

To enable gradient-based optimization through the stochastic sampling operation, the latent vector is sampled using the reparameterization trick~\cite{doersch2016arxiv}.
The latent vector is then sampled according to
\begin{equation}
    \bm{z}
    =
    \bm{\mu}_{\boldsymbol{\phi}}(\bm{x})
    +
    \exp\!\left(\frac12\bm v_\phi(\bm x)\right)
    \odot
    \bm{\epsilon},
    \qquad
    \bm{\epsilon}
    \sim
    \mathcal{N}(\bm{0},\mathbf{I}_d),
    \label{eq:reparameterization}
\end{equation}
where \(\bm{\epsilon}\in\mathbb{R}^{d}\) is a standard Gaussian noise vector and \(\odot\) denotes element-wise multiplication.

The latent input to the drag regressor differs between training and inference.
During training, the regressor receives the latent vector \(\bm{z}\) sampled from the approximate posterior, whereas during inference, the posterior mean is used instead to obtain a deterministic prediction.
Accordingly, the geometric representation \(\mathcal{R}(\mathcal{G})\) introduced in Eq.~\eqref{eq:prediction_formulation}
is defined as
\begin{equation}
    \begin{aligned}
        \mathcal{R}_{\mathrm{train}}(\mathcal{G})
        & = \bm{z},
        & \qquad
        \bm{z}
        & \sim
        q_{\boldsymbol{\phi}}
        \left(
            \bm{z}\mid\mathcal{V}(\mathcal{G})
        \right),
        \\
        \mathcal{R}_{\mathrm{infer}}(\mathcal{G})
        & =
        \bm{\mu}_{\boldsymbol{\phi}}
        \left(
            \mathcal{V}(\mathcal{G})
        \right).
    \end{aligned}
    \label{eq:representation_modes}
\end{equation}

\subsection{Reconstruction and Drag-Prediction Branches}
\label{sec:drag_prediction}

The two downstream branches provide complementary supervision for the latent representation. Reconstruction encourages retention of spatial geometric information, whereas drag regression encourages retention of geometry-dependent information relevant to \(C_{\mathrm d}\).

For reconstruction during training, the decoder \(\mathcal{D}_{\boldsymbol{\theta}}\) maps the sampled latent vector \(\bm{z}\) back to the voxel domain. Its output is a field of pre-sigmoid logits,
\begin{equation}
    \bm{a}
    =
    \mathcal{D}_{\boldsymbol{\theta}}(\bm{z}),
    \label{eq:decoder_logits}
\end{equation}
which can be converted element-wise to reconstructed occupancy probabilities through \(\hat{\bm{x}}=\operatorname{sigmoid}(\bm{a})\).

For drag prediction during training, the regressor combines the sampled latent vector \(\bm z\) with the Reynolds number and evaluates
\begin{equation}
    \hat{C}_{\mathrm{d}}
    =
    f_{\boldsymbol{\psi}}
    \left(
        \bm{z},
        \mathrm{Re}
    \right).
    \label{eq:training_equation}
\end{equation}
Thus, geometry-dependent information is supplied through \(\bm{z}\), while the flow condition enters explicitly through \(\mathrm{Re}\). At inference, the sampled latent vector \(\bm{z}\) is replaced by the posterior mean \(\bm{\mu}_{\boldsymbol{\phi}}\left(\mathcal{V}(\mathcal{G})\right)\) defined in Eq.~\eqref{eq:representation_modes}.

\subsection{Loss Function: Joint Training Objective}
\label{sec:loss_function}
The encoder, decoder, and regressor parameters \((\boldsymbol{\phi},\boldsymbol{\theta},\boldsymbol{\psi})\) are jointly optimized using a total loss comprising a reconstruction term, a KL term, and a drag-prediction term. The joint training objective is:
\begin{equation}
    \mathcal{J}
    =
    w_{\mathrm{rec}}\mathcal{L}_{\mathrm{rec}}
    +
    w_{\mathrm{KL}}\mathcal{L}_{\mathrm{KL}}
    +
    w_{\mathrm{drag}}\mathcal{L}_{\mathrm{drag}},
    \label{eq:joint_training_objective}
\end{equation}
where \(w_{\mathrm{rec}}\), \(w_{\mathrm{KL}}\), and \(w_{\mathrm{drag}}\) control the relative contributions. The three terms play distinct roles in shaping the latent representation. The reconstruction term promotes retention of spatial geometric information, the KL term regularizes the approximate posterior toward the prescribed Gaussian prior, and the drag term promotes retention of variations relevant to the target \(C_{\mathrm{d}}\).
The individual terms are defined below.

Since the input geometry is a sparse binary occupancy field, we formulate the reconstruction objective as a positive-class-weighted binary cross-entropy loss evaluated directly on the decoder's pre-sigmoid logits. For voxel \(j\) of sample \(i\), this contribution is
\begin{equation}
    \ell_{\mathrm{rec}}(a_{ij},x_{ij})
    =
    w_{\mathrm{pos}}x_{ij}\log\left(1+\exp(-a_{ij})\right)
    +
    (1-x_{ij})\log\left(1+\exp(a_{ij})\right),
\end{equation}
where \(a_{ij}\) is the pre-sigmoid logit, \(x_{ij}\in\{0,1\}\) is the ground-truth occupancy, and \(w_{\mathrm{pos}}\) is the positive-class weight. Because occupied voxels constitute only a small fraction of the grid, \(w_{\mathrm{pos}}>1\) assigns greater weight to errors on occupied voxels and reduces the influence of the empty background on the reconstruction objective~\cite{brock2016arxiv}.
The mean reconstruction loss is then
\begin{equation}
    \mathcal{L}_{\mathrm{rec}}
    =
    \frac{1}{B}
    \sum_{i=1}^{B}
    \frac{1}{N_{\mathrm{vox}}}
    \sum_{j=1}^{N_{\mathrm{vox}}}
    \ell_{\mathrm{rec}}
    \left(
        a_{ij},x_{ij}
    \right),
    \label{eq:reconstruction_loss}
\end{equation}
where \(B\) is the batch size and \(N_{\mathrm{vox}}=N_v^3\) is the total number of voxels per sample.

Following the standard VAE formulation, we regularize the approximate posterior distribution toward a standard Gaussian prior \(p(\bm{z})=\mathcal{N}(\bm{0},\mathbf{I}_d)\). Because the approximate posterior is modeled as a diagonal Gaussian, the KL divergence can be decomposed across dimensions and evaluated in closed form. Let \(\mu_{ik}=[\bm{\mu}_{\boldsymbol{\phi}}(\bm{x}_i)]_k\) and \(v_{ik}=[\bm{v}_{\boldsymbol{\phi}}(\bm{x}_i)]_k\) denote the mean and log-variance of latent dimension \(k\) for sample \(i\). The KL contribution for this dimension is
\begin{equation}
    \ell_{\mathrm{KL},ik}
    =
    \frac{1}{2}
    \left[
        \mu_{ik}^{2}
        +
        \exp(v_{ik})
        -
        1
        -
        v_{ik}
    \right].
\end{equation}
Therefore, the KL term is defined as the empirical mean over the mini-batch and latent dimensions:
\begin{equation}
    \mathcal{L}_{\mathrm{KL}}
    =
    \frac{1}{B}
    \sum_{i=1}^{B}
    \frac{1}{d}
    \sum_{k=1}^{d}
    \ell_{\mathrm{KL},ik}.
    \label{eq:kl_loss}
\end{equation}
This normalization reduces the direct dependence of the numerical scale of the KL term on the mini-batch size and latent dimensionality \(d\).

The supervised drag-prediction loss is formulated as a target-magnitude-weighted squared error:
\begin{equation}
    \mathcal{L}_{\mathrm{drag}}
    =
    \frac{1}{B}
    \sum_{i=1}^{B}
    \frac{
        \left[
            C_{\mathrm{d},i}
            -
            f_{\boldsymbol{\psi}}
            \left(
                \bm{z}_i,
                \mathrm{Re}_i
            \right)
        \right]^2
    }
    {\left(|C_{\mathrm{d},i}|+\delta\right)^{\alpha}},
    \label{eq:drag_loss}
\end{equation}
where \(\bm{z}_i\) and \(\mathrm{Re}_i\) are the latent sample and Reynolds number of sample \(i\), respectively, and \(C_{\mathrm{d},i}\) is the corresponding target drag coefficient. The small constant \(\delta=10^{-6}\) avoids a zero denominator, and \(\alpha\) controls the dependence of the weight on the target magnitude. Setting \(\alpha=0\) recovers the unweighted mean squared error, whereas \(\alpha=2\) approximates relative squared-error weighting. The numerical values of the loss weights and other training hyperparameters are specified in Section~\ref{sec:network_architecture_training}.

The reparameterization enables gradients from both downstream objectives to propagate back to the encoder through the sampled latent vector \(\bm z\). Specifically, \(\mathcal{L}_{\mathrm{rec}}\) backpropagates through the decoder, whereas \(\mathcal{L}_{\mathrm{drag}}\) backpropagates through the drag regressor. \(\mathcal{L}_{\mathrm{KL}}\) acts directly on the approximate posterior distribution parameterized by the encoder.

\section{Experimental Setup}

We systematically evaluate our methodology with a group of experiments, E1, E2a, E2b, and E3. These experiments compare drag predictors across increasing geometric complexity and different training coverage of geometric variations.

\subsection{Particle geometry generation and voxelization}
\label{sec:geometry_generation}
\subsubsection{Particle geometry generation}

Particle geometries are generated using an ellipsoid-specific adaptation of the open-source implementation of spherical-harmonic (SH) particle-generation method of Wei et al.~\cite{wei2018particle}.
The base ellipsoid is parameterized by the aspect ratio \(\mathrm{Ar}\), defined as the longitudinal diameter divided by the transverse diameter, and is normalized to a nominal volume. The sampled range \(\mathrm{Ar}\in[0.25,2.5]\) spans oblate (\(\mathrm{Ar}<1\)), spherical (\(\mathrm{Ar}=1\)), and prolate (\(\mathrm{Ar}>1\)) shapes.

\begin{figure}[htbp]
    \centering
    \includegraphics[width=0.72\linewidth]{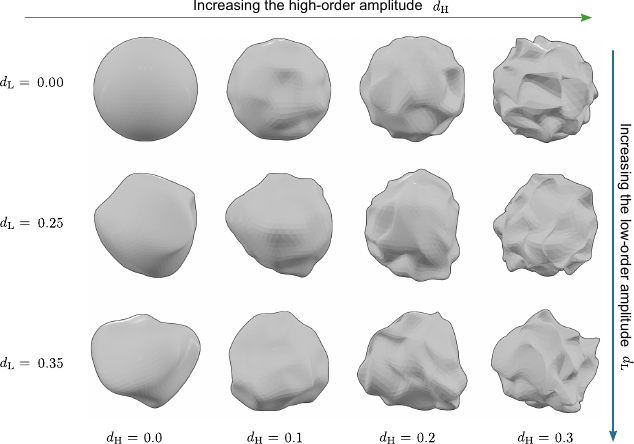}
    \caption{An example of a spherical base shape (\(\mathrm{Ar}=1\)) perturbed by different pairs of \(d_{\mathrm{L}}\) and \(d_{\mathrm{H}}\).}
    \label{fig:controlled_particle_shapes}
\end{figure}

Geometric irregularity is introduced by perturbing the spherical-harmonic (SH) representation of a base ellipsoid. The dimensionless parameter \(d_{\mathrm L}\) controls perturbations in degrees \(n=2,\ldots,8\). These perturbations mainly change the overall particle shape. The dimensionless parameter \(d_{\mathrm H}\) controls perturbations in degrees \(n=9,\ldots,14\). These perturbations add finer surface features. The generator uses the Frobenius norm of the corresponding base ellipsoid’s SH coefficient matrix as the reference scale for both parameters. Figure~\ref{fig:controlled_particle_shapes} illustrates their effects for a spherical base shape. Appendix~\ref{app:geometry_generation_details} gives the coefficient-scaling procedure.

For each geometry, the SH coefficients are independently sampled according to the prescribed \(d_{\mathrm{L}}\) and \(d_{\mathrm{H}}\) conditions. Geometries generated at different perturbation levels therefore correspond to independent realizations rather than a nested sequence of shared SH coefficients. The SH perturbation is applied directly to the nominally volume-normalized base ellipsoid, without an additional volume-rescaling step, so the resulting particle volume may differ slightly from the nominal value. Each perturbed geometry is then rotated about the \(y\)-axis to impose an incident angle \(\alpha_{\mathrm{inc}}\in[0^\circ,90^\circ]\), with \(\alpha_{\mathrm{inc}}=0^\circ\) corresponding to alignment of the particle \(x\)-axis with the free-stream direction. The parameters \(\mathrm{Ar}\), \(d_{\mathrm{L}}\), and \(d_{\mathrm{H}}\) refer to the geometry before this orientation is applied.

\subsubsection{Voxel occupancy fields}
\label{sec:voxelization}
Voxelization converts the oriented particle mesh into a binary voxel representation for input to the three-dimensional convolutional network.
Prior to voxelization, each oriented mesh is translated so that the center of its axis-aligned bounding box coincides with the origin.
The mesh is then scaled by \(s=0.9/L_{\max}\), where \(L_{\max}\) is the largest bounding-box dimension, such that its longest dimension becomes 0.9 in the normalized coordinate system. Surface voxels are generated using a subdivision-based voxelization procedure, after which the enclosed interior is filled to obtain a binary occupancy field, with occupied and empty voxels assigned values of 1 and 0, respectively. The resulting voxel matrix is centered within a fixed \(64\times64\times64\) array, with the remaining entries zero-padded, yielding a single-channel network input of size \(1\times64\times64\times64\).

\subsection{Network architecture and training}
\label{sec:network_architecture_training}

The VAE-MLP formulation described in Section~\ref{sec:overall_framework} is implemented using a three-dimensional convolutional encoder and decoder together with an MLP drag regressor. All main experiments use the same overall architecture; the only architectural variable is the latent dimensionality. E1 uses \(d=2\) to permit direct visualization of the latent space, whereas E2a, E2b, and E3 use \(d=16\), selected based on the latent-dimension sensitivity analysis in~\ref{app:latent_sensitivity}. The complete layer configuration and tensor shapes throughout the network are provided in~\ref{app:architecture_specification}.

The loss configuration is held constant across all VAE-MLP experiments. The reconstruction, KL, and drag terms are weighted by \((w_{\mathrm{rec}},w_{\mathrm{KL}},w_{\mathrm{drag}})=(1.0,0.1,10.0)\), with these weights fixed based on preliminary training trials. The positive-voxel weight \(w_{\mathrm{pos}}=17\) compensates for the imbalance between occupied and empty voxels. The drag-loss exponent \(\alpha=0.4\) reduces the dominance of larger \(C_{\mathrm d}\) values in the training objective.

The models are trained for 200 epochs using Adam with a learning rate of \(10^{-3}\) and a batch size of 8. The checkpoint with the minimum validation drag loss is retained for external evaluation. Both models use seed 42 in E1, E2a, and E2b, while E3 is repeated with seeds 42--46 using a fixed train--validation partition.
Additional details on the software environment and random-seed configuration are provided in~\ref{app:training_reproducibility}.

\subsection{CFD reference datasets}
\label{sec:cfd_reference_data}

The drag coefficients used as learning targets are obtained from transient incompressible-flow simulations performed in OpenFOAM v2306~\cite{weller1998tensorial}.
The simulations cover the Reynolds-number range \(\mathrm{Re}\in[10,200]\).
For each particle geometry, \(\mathrm{Re}\) is defined using its geometry-specific volume-equivalent diameter \(D_{\mathrm{eq}}\), with the prescribed value controlled by varying the fluid kinematic viscosity.
For each particle, the drag-reference area is the cross-sectional area
of its volume-equivalent sphere,
\begin{equation}
    A_{\mathrm{ref}}=\frac{\pi D_{\mathrm{eq}}^2}{4}.
    \label{eq:reference_area}
\end{equation}
Here, \(D_{\mathrm{eq}}\) is determined from the volume of the final
particle geometry. The reference area therefore depends on particle
volume and is invariant under changes in particle orientation.
During the simulations, the hydrodynamic force acting on each particle is evaluated by integrating the pressure and skin-friction contributions over the particle-surface patches~\cite{openfoam_forceCoeffs}.
The scalar learning target for each configuration is obtained by taking the arithmetic mean over the final portion of the transient drag coefficient history.
The complete mathematical formulation, computational-domain dimensions, boundary conditions, numerical schemes, and drag-evaluation procedure are demonstrated in~\ref{app:cfd_theory} and~\ref{app:cfd_drag_evaluation}.

The numerical framework and mesh resolution are assessed through a three-level mesh-verification study of regular particles, with the full procedure and quantitative results reported in~\ref{app:cfd_mesh_verification}. The selected medium-grid shows good agreement with the fine-grid solutions, with a mean relative difference of about \(0.37\%\) across representative combinations of particle shape, orientation, and Reynolds number, while the mesh-related differences (\(\delta_{\mathrm{M/F}}\)) remain below \(0.64\%\). These comparisons support the numerical consistency of the selected production grid for the regular-particle cases examined.

\subsection{Dataset partition and experimental regimes}

Each training dataset is constructed by first sampling intrinsic particle geometries and then associating each geometry with multiple flow conditions. Training geometries are generated using Latin hypercube sampling (LHS) within the parameter domain prescribed for each experiment, while the corresponding incident angles and Reynolds numbers are jointly sampled using a two-dimensional LHS. The partition between train and validation is performed at the intrinsic-geometry level, such that all orientations and flow conditions derived from the same intrinsic geometry are assigned to the same subset. This geometry-disjoint partition prevents geometry-level leakage between training and validation. All reported prediction-performance metrics are evaluated on independently generated external datasets. The exact training domains, sample counts, train- and validation- configurations are specified in~\ref{app:dataset_configuration}.

\begin{figure}[htbp]
    \centering
    \includegraphics[width=0.52\linewidth]{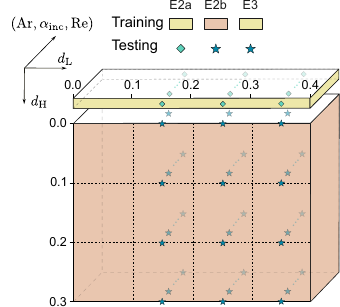}
    \caption{Training and external-evaluation distributions for E2a, E2b, and E3 in the \((d_{\mathrm{L}},d_{\mathrm{H}})\) parameter space.}
    \label{fig:dataset_schematic}
\end{figure}

E1 uses regular ellipsoids, and E2a uses particles with low-order perturbations. E2b and E3 differ in their training coverage of high-order perturbations: E2b includes \(d_{\mathrm H}\in[0,0.3]\), whereas E3 is trained with \(d_{\mathrm H}=0\). Figure~\ref{fig:dataset_schematic} illustrates the relationships between the training and evaluation distributions.

The external evaluation datasets are constructed from prescribed Cartesian-product grids to enable controlled comparisons across geometric and flow conditions. E2b and E3 share a complete external evaluation grid comprising 864 cases. Aggregate prediction metrics use all 864 cases in E2b and the 648 cases with \(d_{\mathrm H}>0\) in E3. The remaining 216 cases with \(d_{\mathrm H}=0\) serve as condition-matched references for the E3 analyses of drag differences and changes in geometric representations. Table~\ref{tab:external_config} lists the grids used for aggregate predictive evaluation.

\begin{table}[htbp]
    \centering
    \caption{Prescribed Cartesian-product grids for external evaluation. The E3 row lists the 648 cases used for aggregate predictive evaluation; an additional 216 cases at \(d_{\mathrm H}=0\), with the same values of the other parameters, serve as condition-matched references.  Each sample corresponds to one combination of geometry parameters, incident angle, and Reynolds number. \(S_{\mathrm{Ar}}=\{0.50,0.75,1.50,2.00\}\). The Reynolds-number sets are \(S_{\mathrm{Re}}=\{12.5,15,22.5,30,45,75,100,150,200\}\).}
    \label{tab:external_config}
    \small
    \begin{tabular}{lcccccc}
        \toprule
        \textbf{Exp.}
        & \(\mathrm{Ar}\)
        & \(d_{\mathrm{L}}\)
        & \(d_{\mathrm{H}}\)
        & \(\alpha_{\mathrm{inc}}\ (^\circ)\)
        & \(\mathrm{Re}\)
        & Total samples \\
        \midrule
        E1
        & \(S_{\mathrm{Ar}}\)
        & 0.00
        & 0.0
        & \(\{0,30,45,60,90\}\)
        & \(S_\mathrm{Re}\)
        & 180 \\
        \addlinespace
        E2a
        & \(S_{\mathrm{Ar}}\)
        & \(\{0.15,0.25,0.35\}\)
        & 0.0
        & \(\{0,30,60\}\)
        & \(S_\mathrm{Re}\)
        & 324 \\
        \addlinespace
        E2b
        & \(S_{\mathrm{Ar}}\)
        & \(\{0.25,0.35\}\)
        & \(\{0.0,0.1,0.2,0.3\}\)
        & \(\{0,30,60\}\)
        & \(S_\mathrm{Re}\)
        & 864 \\
        \addlinespace
        E3
        & \(S_{\mathrm{Ar}}\)
        & \(\{0.25,0.35\}\)
        & \(\{0.1,0.2,0.3\}\)
        & \(\{0,30,60\}\)
        & \(S_\mathrm{Re}\)
        & 648 \\
        \bottomrule
    \end{tabular}
\end{table}

\subsection{Comparison baselines}
\label{sec:comparison_baselines}

Two types of baselines are used: Descriptor-MLP and empirical drag correlations. Descriptor-MLP provides a comparison between hand-crafted and learned geometric representations within a neural-network predictor. The Ke and H\"{o}lzer--Sommerfeld correlations provide established empirical references for regular ellipsoids and irregular particles, respectively.

Descriptor-MLP takes the volume-equivalent diameter \(D_{\mathrm{eq}}\), sphericity \(\Phi\), crosswise sphericity \(\Phi_{\mathrm{cross}}\), and Reynolds number as inputs. The three geometric descriptors characterize particle size, surface compactness, and orientation-dependent projected geometry~\cite{holzer2008}, and are defined as
\begin{equation}
    D_{\mathrm{eq}}
    =
    \left(\frac{6V}{\pi}\right)^{1/3},
    \qquad
    \Phi
    =
    \frac{\pi D_{\mathrm{eq}}^2}{A_p},
    \qquad
    \Phi_{\mathrm{cross}}
    =
    \frac{\pi D_{\mathrm{eq}}^2}{4A_\perp},
    \label{eq:descriptor_definitions}
\end{equation}
where \(V\) and \(A_p\) denote the particle volume and surface area, respectively, and \(A_\perp\) is the projected silhouette area normal to the incident-flow direction.
The surface area is computed by summing the triangular facet areas, while the projected area is computed from the union of the orthogonally projected facets, counting overlapping regions only once.

Descriptor-MLP uses the same regression backbone as the drag-prediction branch of VAE-MLP, with the latent input replaced by the geometric descriptors.
Both predictors use the same drag-prediction loss, optimizer settings, number of training epochs, and checkpoint-selection criterion.
Their architecture details are provided in~\ref{app:architecture_specification}.

The ellipsoid-specific Ke correlation~\cite{ke2018} is used in E1, and the crosswise-sphericity form of the H\"{o}lzer--Sommerfeld correlation~\cite{holzer2008} is used in E2a and E2b.
Both correlations are evaluated using their published coefficients without recalibration to the present CFD data. Their formulations and coefficients are given in~\ref{app:empirical_baseline_formulations}.

\subsection{Evaluation metrics}

We quantify prediction error using the absolute relative error \(\varepsilon_{\mathrm{rel},i}\). The mean absolute relative error (MARE) is the primary aggregate metric, while \(\max_i\varepsilon_{\mathrm{rel},i}\) characterizes the largest individual deviation. These quantities are defined as
\begin{equation}
    e_{\mathrm{rel},i}
    =
    \frac{
        \hat{C}_{\mathrm{d},i}-C_{\mathrm{d},i}
    }{
        C_{\mathrm{d},i}
    },
    \qquad
    \varepsilon_{\mathrm{rel},i}
    =
    \frac{
        \left|\hat{C}_{\mathrm{d},i}-C_{\mathrm{d},i}\right|
    }{
        C_{\mathrm{d},i}
    },
    \qquad
    \mathrm{MARE}
    =
    \frac{1}{N_s}
    \sum_{i=1}^{N_s}
    \varepsilon_{\mathrm{rel},i},
\end{equation}
where \(N_s\) denotes the number of samples included in the evaluation. All relative errors are reported as percentages.

For condition-resolved analyses, we define the condition-wise MARE (cMARE) as the MARE computed over samples restricted to a fixed parameter set \((\mathrm{Re},d_{\mathrm{L}},d_{\mathrm{H}})\), rather than over the entire evaluation dataset. The coefficient of determination \(R^2\) is reported as a complementary measure of overall predictive agreement. All cross-model comparisons use identical evaluation samples and configurations to ensure a common basis for comparison. Additional analysis-specific metrics are defined locally in the sections where they are introduced.

\subsection{Experimental matrix}
The geometry regimes progress from regular ellipsoids in E1 to low-order irregular particles in E2a and then to both high- and low-order irregular particles in E2b, which contain both low- and high-order perturbations (\(d_{\mathrm L}\) and \(d_{\mathrm H}\)). E3 retains low-order irregularity in the training geometries but excludes high-order perturbations (\(d_{\mathrm H}=0\)) from training, while evaluating particles with \(d_{\mathrm H}>0\). Table~\ref{tab:exp_matrix} summarizes the training and evaluation regimes, experimental purposes, and compared methods for the four main experiments.

\begin{table}[htbp]
    \centering
    \small
    \caption{Experimental matrix for the main predictive evaluation. Here, \(d_{\mathrm {L,H}}=0\) indicates the absence of the corresponding perturbation in the dataset, whereas \(d_{\mathrm {L,H}}>0\) indicates its presence in the dataset. H--S denotes the Hölzer--Sommerfeld correlation; corr. denotes correlation.}
    \label{tab:exp_matrix}
    \begin{tabularx}{\textwidth}{
            >{\centering\arraybackslash}m{0.055\textwidth}
            >{\centering\arraybackslash}m{0.16\textwidth}
            >{\centering\arraybackslash}m{0.16\textwidth}
            >{\centering\arraybackslash}m{0.27\textwidth}
            >{\centering\arraybackslash}m{0.25\textwidth}
        }
        \toprule
        \textbf{No.}                              &
        \textbf{Training regime}                  &
        \textbf{Evaluation regime}                &
        \textbf{Purpose}                          &
        \textbf{Compared methods}                   \\
        \midrule

        E1                                        &
        \(d_{\mathrm L}=0,\ d_{\mathrm H}=0\)     &
        \(d_{\mathrm L}=0,\ d_{\mathrm H}=0\)     &
        Validation on ellipsoids                  &
        \shortstack[c]{VAE-MLP,                     \\ Descriptor-MLP,\\ Ke et al. corr.} \\

        \addlinespace[3pt]

        E2a                                       &
        \(d_{\mathrm L}>0,\ d_{\mathrm H}=0\)     &
        \(d_{\mathrm L}>0,\ d_{\mathrm H}=0\)     &
        Low-order perturbations covered in training           &
        \shortstack{VAE-MLP,                        \\ Descriptor-MLP,\\ H--S. corr.} \\

        \addlinespace[3pt]

        E2b                                       &
        \(d_{\mathrm L}>0,\ d_{\mathrm H}>0\)     &
        \(d_{\mathrm L}>0,\ d_{\mathrm H}\geq 0\) &
        Low and high-order perturbations covered in training           &
        \shortstack{VAE-MLP,                        \\ Descriptor-MLP,\\ H--S. corr.} \\

        \addlinespace[3pt]

        E3                                        &
        \(d_{\mathrm L}>0,\ d_{\mathrm H}=0\)     &
        \(d_{\mathrm L}>0,\ d_{\mathrm H}>0\)     &
        High-order perturbations absent from training                             &
        \shortstack{VAE-MLP,                        \\ Descriptor-MLP} \\

        \bottomrule
    \end{tabularx}
\end{table}

\section{Results and Discussion}
We assess whether a latent representation learned from three-dimensional particle geometry improves drag prediction and generalization to unseen geometric changes compared with hand-crafted descriptors. E1 establishes a regular-ellipsoid baseline for comparing drag-prediction accuracy. E2a and E2b compare the models in two irregular-geometry regimes whose evaluated variations are represented during training: low-order perturbations alone and combined low- and high-order perturbations. E3 evaluates drag-prediction generalization to the tested high-order surface perturbations excluded from training. The corresponding CFD drag variation across the controlled perturbation space is reported in~\ref{app:drag_variation}.

\subsection{E1: Learned latent representation and drag coefficient prediction for regular ellipsoids}

E1 evaluates the learned latent representation and drag prediction on regular ellipsoids. We first investigate how aspect ratios and incident angles, which govern ellipsoid shapes and orientations, are represented in the latent space. We then evaluate the drag predictions from VAE-MLP, Descriptor-MLP, and the Ke correlation~\cite{ke2018}.

\begin{figure}[htbp]
    \centering
    \includegraphics[width=1.0\linewidth]{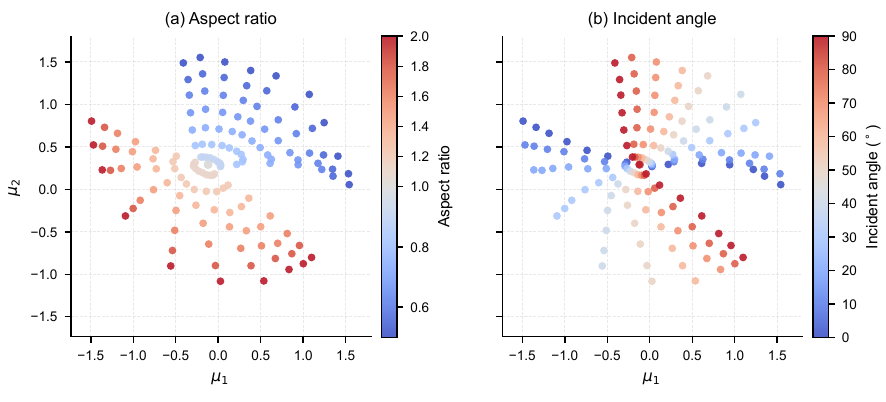}
    \caption{Posterior means of the jointly trained two-dimensional E1 latent space, colored by (a) aspect ratio and (b) incident angle.}
    \label{fig:latent-space_e1}
\end{figure}

We visualize the latent space using a separate geometry-only grid spanning aspect ratio and incident angle. This grid has no CFD drag targets. Figure~\ref{fig:latent-space_e1} shows the posterior means \(\mu_{1}\) and \(\mu_2\) produced by the encoder for these geometries. In Figure~\ref{fig:latent-space_e1}(a), near-spherical geometries (\(\mathrm{Ar}\approx1\)) cluster in the central region of the latent distribution, whereas increasingly oblate and prolate geometries extend along distinct trajectories. In Figure~\ref{fig:latent-space_e1}(b), the progression from blue to red along the trajectories corresponds to increasing incident angle from \(0^\circ\) to \(90^\circ\), showing that the latent position also varies continuously with orientation. This organization indicates that the two-dimensional latent representation retains information associated with both controlled geometric variables, although neither is supplied as an explicit scalar input.

\begin{figure}[htbp]
    \centering
    \includegraphics[width=1\linewidth]{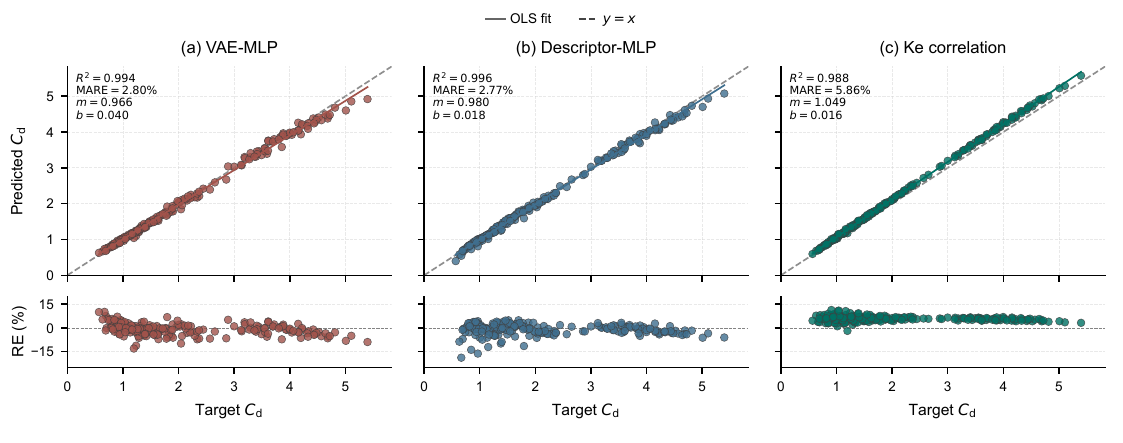}
    \caption{Drag predictions and signed relative errors on the E1 external evaluation set for (a) VAE-MLP, (b) Descriptor-MLP, and (c) the Ke correlation~\cite{ke2018}. Solid colored lines denote ordinary least-squares fits, and gray dashed lines denote ideal parity.}
    \label{fig:parity_e1}
\end{figure}

VAE-MLP and Descriptor-MLP achieve similar average predictive accuracy on the E1 external evaluation set. Their MAREs, listed in Table~\ref{tab:e1_prediction_performance}, are \(2.80\%\) and \(2.77\%\), respectively. Their maximum absolute relative errors differ more substantially: \(12.97\%\) for VAE-MLP and \(30.37\%\) for Descriptor-MLP.
The signed relative errors in Figure~\ref{fig:parity_e1} show greater dispersion for Descriptor-MLP at \(C_{\mathrm d}<1.5\), consistent with its larger worst-case error.
The established Ke correlation gives a higher MARE of \(5.86\%\), while its maximum absolute relative error is \(11.37\%\). Its predominantly positive signed relative errors indicate a systematic tendency to overpredict \(C_{\mathrm d}\) on the present evaluation set.

\begin{table}[htbp]
    \centering
    \small
    \caption{Global prediction metrics on the E1 external evaluation set. Max. \(\varepsilon_{\mathrm{rel}}\) denotes the maximum absolute relative error.}
    \label{tab:e1_prediction_performance}
    \begin{tabular}{lccc}
        \toprule
        Model          & \(R^2\) & MARE (\%) & Max. \(\varepsilon_{\mathrm{rel}}\) (\%) \\
        \midrule
        VAE-MLP        & 0.994   & 2.80      & 12.97 \\
        Descriptor-MLP & 0.996   & 2.77      & 30.37 \\
        Ke correlation & 0.988   & 5.86      & 11.37 \\
        \bottomrule
    \end{tabular}
\end{table}

Together, the latent-space patterns and comparable average prediction errors on regular ellipsoids provide a baseline for the irregular-geometry comparisons that follow.
\subsection{E2a and E2b: Interpolative study under geometric perturbations}

E2a and E2b experiments assess predictive performance of models in two geometric regimes represented during training. E2a includes low-order perturbations alone, whereas E2b additionally includes high-order surface perturbations.

In E2a, VAE-MLP achieves a lower MARE than Descriptor-MLP, \(3.37\%\) compared with \(4.39\%\). Table~\ref{tab:e2a_prediction_performance} also shows a lower maximum absolute relative error for VAE-MLP: \(14.32\%\) compared with \(24.83\%\) for Descriptor-MLP. Both learned models also give lower errors than the H\"{o}lzer--Sommerfeld correlation, whose MARE is \(9.18\%\).

\begin{table}[htbp]
    \centering
    \small
    \caption{Global external-evaluation prediction metrics in E2a. Max. \(\varepsilon_{\mathrm{rel}}\) denotes the maximum absolute relative error.}
    \label{tab:e2a_prediction_performance}
    \begin{tabular}{lccc}
        \toprule
        Model                  & \(R^2\)        & MARE (\%)     & Max. \(\varepsilon_{\mathrm{rel}}\) (\%) \\
        \midrule
        VAE-MLP                & \textbf{0.993} & \textbf{3.37} & \textbf{14.32} \\
        Descriptor-MLP         & 0.992          & 4.39          & 24.83 \\
        H\"{o}lzer--Sommerfeld & 0.969          & 9.18          & 38.14 \\
        \bottomrule
    \end{tabular}
\end{table}

\begin{figure}[htbp]
    \centering
    \includegraphics[width=1.0\linewidth]{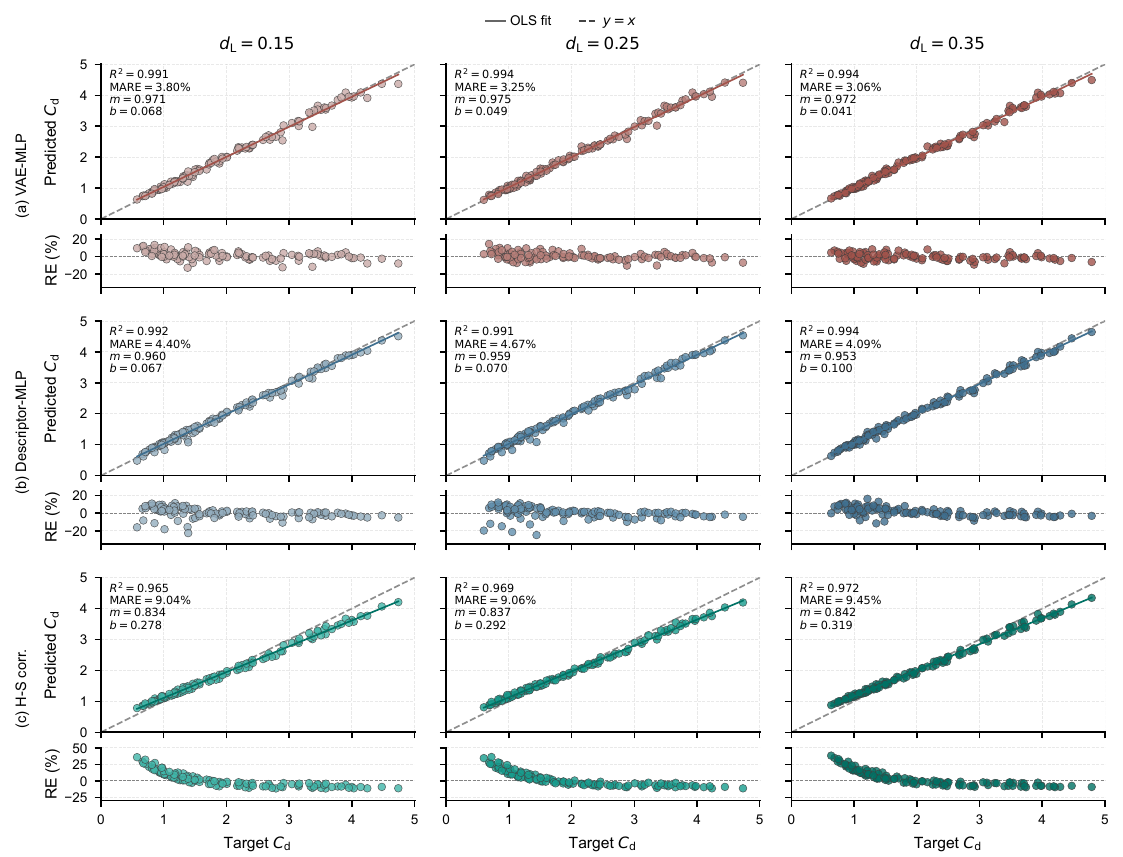}
    \caption{E2a drag predictions and signed relative errors resolved by \(d_{\mathrm L}\) level. Columns correspond to \(d_{\mathrm{L}}=0.15\), \(0.25\), and \(0.35\), while rows correspond to (a) VAE-MLP, (b) Descriptor-MLP, and (c) the H\"{o}lzer--Sommerfeld correlation. Solid lines denote panel-specific ordinary least-squares fits, and gray dashed lines denote ideal parity. The lower panels show the signed relative error (RE).}
    \label{fig:e2a_parity}
\end{figure}

Figure~\ref{fig:e2a_parity} shows that VAE-MLP has a lower MARE than Descriptor-MLP at \(d_{\mathrm L}=0.15\), \(0.25\), and \(0.35\). The H\"{o}lzer--Sommerfeld correlation exhibits a consistent drag-dependent bias across these levels: its fitted slopes are below unity, and its signed relative errors change from positive at lower \(C_{\mathrm d}\) to negative at higher \(C_{\mathrm d}\).

\begin{figure}[htbp]
    \centering
    \includegraphics[width=1.0\linewidth]{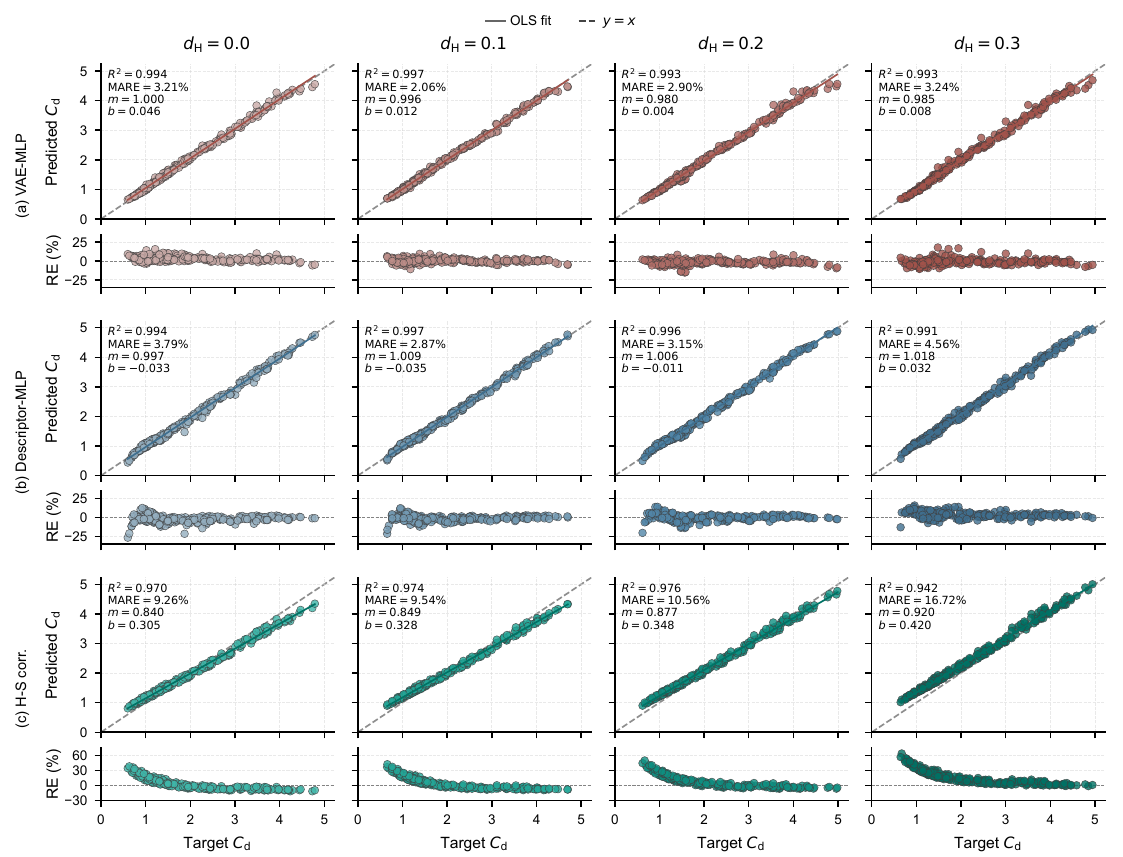}
    \caption{E2b experiment drag predictions and signed relative errors for geometries whose low- and high-order perturbations are both represented in training.}
    \label{fig:e2b_parity}
\end{figure}

E2b extends this comparison to particles containing both low- and high-order perturbations. Table~\ref{tab:e2b_prediction_performance} presents the aggregate E2b results. VAE-MLP and Descriptor-MLP attain the same \(R^2=0.994\), with MAREs of \(2.85\%\) and \(3.59\%\), respectively. VAE-MLP also has a lower maximum absolute relative error (\(17.85\%\) versus \(26.60\%\)). Figure~\ref{fig:e2b_parity} shows a lower VAE-MLP MARE at each evaluated \(d_{\mathrm H}\) level, with the largest difference between the models at \(d_{\mathrm H}=0.3\). The H\"{o}lzer--Sommerfeld correlation exhibits substantially larger errors.

\begin{table}[htbp]
    \centering
    \footnotesize
    \setlength{\tabcolsep}{4pt}
    \caption{External-evaluation prediction metrics in E2b. Max. \(\varepsilon_{\mathrm{rel}}\) denotes the maximum absolute relative error. cMARE is evaluated over samples sharing fixed \((\mathrm{Re},d_{\mathrm{L}},d_{\mathrm{H}})\) values, and its maximum across the tested conditions is reported.}
    \label{tab:e2b_prediction_performance}
    \begin{tabular}{lcccc}
        \toprule
        Model                  & \(R^2\)        & MARE (\%)     & Max. \(\varepsilon_{\mathrm{rel}}\) (\%) & Maximum cMARE (\%) \\
        \midrule
        VAE-MLP                & \textbf{0.994} & \textbf{2.85} & \textbf{17.85}                           & \textbf{5.73} \\
        Descriptor-MLP         & \textbf{0.994} & 3.59          & 26.60                                    & 7.28 \\
        H\"{o}lzer--Sommerfeld & 0.966          & 11.52         & 63.74                                    & 39.46 \\
        \bottomrule
    \end{tabular}
\end{table}

\begin{figure}[htbp]
    \centering
    \includegraphics[width=1.0\linewidth]{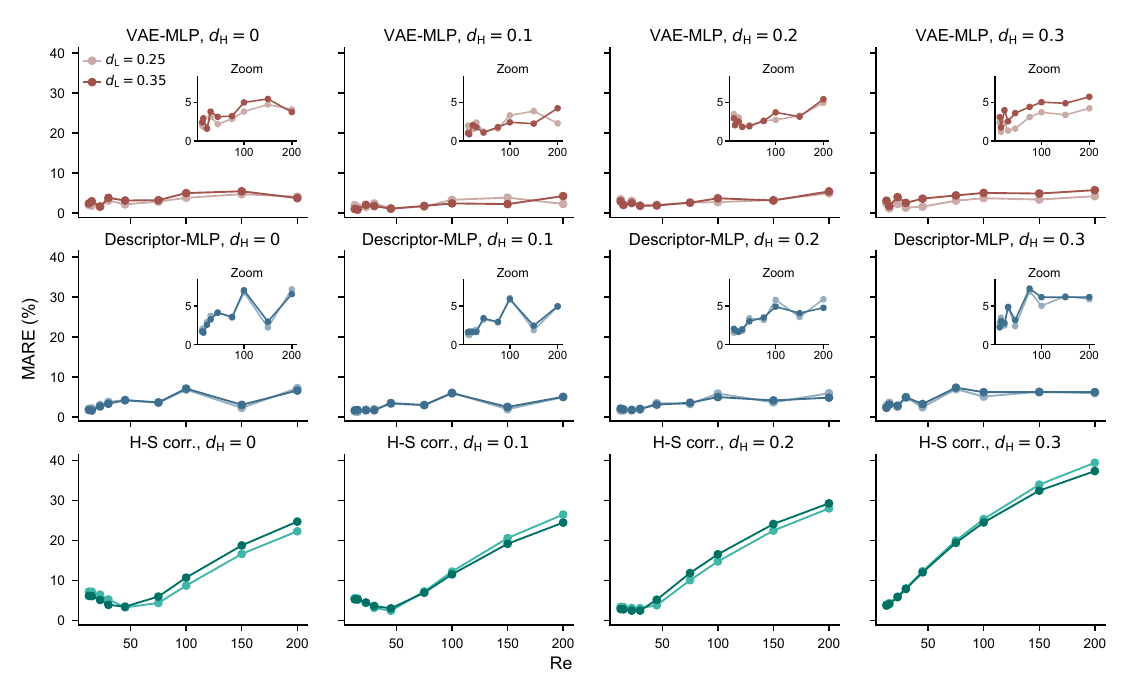}
    \caption{Reynolds-number-resolved MARE in E2b, grouped by \(d_{\mathrm{H}}\) for VAE-MLP, Descriptor-MLP, and H\"{o}lzer--Sommerfeld (HS).}
    \label{fig:e2b_mare_vs_re_by_dh}
\end{figure}

Figure~\ref{fig:e2b_mare_vs_re_by_dh} shows that VAE-MLP has a modestly narrower condition-level error range across the evaluated Reynolds numbers and high-order perturbation levels. Its maximum cMARE is \(5.73\%\), compared with \(7.28\%\) for Descriptor-MLP. Complementary results grouped by \(d_{\mathrm L}\), together with condition-wise error heatmaps over the full \((d_{\mathrm L},d_{\mathrm H})\) evaluation grid, are provided in~\ref{app:e2b_additional_condition_results}.

VAE-MLP yields lower aggregate errors than Descriptor-MLP in both E2a and E2b, where the evaluated geometric variations are represented during training. These results suggest that learning features from three-dimensional geometry is useful for drag prediction within the two evaluated irregular-geometry regimes. E3 tests whether this advantage extends to high-order perturbations absent from training.

\subsection{E3: Extrapolative study under unseen high-order perturbations}
\label{sec:e3_extrapolative_study}

E3 evaluates both models on high-order surface perturbations absent from training. We first compare the prediction errors of VAE-MLP and Descriptor-MLP across different perturbation levels and flow conditions. We then examine condition-matched drag differences and conduct input-substitution analyses to investigate the prediction behavior of the two models.

\subsubsection{Predictive performance}
To assess the stability of the E3 model comparison, we train both models with seeds 42--46 while holding the data split fixed.
VAE-MLP achieves lower aggregate MARE than Descriptor-MLP in all five E3 training runs, as shown in Table~\ref{tab:e3_multiseed_performance}.
Aggregate prediction metrics are calculated over the 648 cases with \(d_{\mathrm H}=0.1\), \(0.2\), and \(0.3\), evaluated at \(d_{\mathrm L}=0.25\) and \(0.35\); the 216 cases at \(d_{\mathrm H}=0\) serve as matched references for the response analyses.
The cross-seed mean MARE is \(3.99\%\pm0.55\%\) for VAE-MLP and \(6.11\%\pm1.72\%\) for Descriptor-MLP.
At \(d_{\mathrm H}=0.3\), VAE-MLP also has lower MARE in every run, with corresponding cross-seed values of \(4.97\%\pm0.78\%\) and \(9.88\%\pm3.22\%\).
Values are expressed as the mean \(\pm\) sample standard deviation across training seeds.

\begin{table}[htbp]
    \centering
    \small
    \setlength{\tabcolsep}{4pt}
    \caption{E3 MARE (\%) across training seeds with the data split fixed at seed 42. Global metrics are calculated over the 648 cases with \(d_{\mathrm H}>0\). The final row gives the cross-seed arithmetic mean and sample standard deviation.}
    \label{tab:e3_multiseed_performance}
    \begin{tabular}{lcccc}
        \toprule
                        & \multicolumn{2}{c}{Global} & \multicolumn{2}{c}{\(d_{\mathrm H}=0.3\)} \\
        \cmidrule(lr){2-3}\cmidrule(lr){4-5}
        Training seed   & VAE-MLP                    & Descriptor-MLP & VAE-MLP & Descriptor-MLP \\
        \midrule
        42              & 4.545                      & 8.865          & 5.969   & 15.215 \\
        43              & 3.508                      & 4.729          & 4.337   & 7.811 \\
        44              & 4.602                      & 6.701          & 5.493   & 10.664 \\
        45              & 3.799                      & 5.100          & 4.920   & 7.758 \\
        46              & 3.490                      & 5.161          & 4.108   & 7.969 \\
        \midrule
        Mean \(\pm\) SD & \(3.989\pm0.548\)          & \(6.111\pm1.715\)
                        & \(4.966\pm0.778\)          & \(9.884\pm3.221\) \\
        \bottomrule
    \end{tabular}
\end{table}

Among the five runs listed in Table~\ref{tab:e3_multiseed_performance}, seed 46 gives the median MARE difference between the models over the 648 unseen-perturbation cases and at \(d_{\mathrm H}=0.3\). We use this run to illustrate condition-resolved error patterns in the figures below.
Figure~\ref{fig:e3_parity} shows that VAE-MLP has lower MARE at all three unseen perturbation levels in this run.
The largest difference occurs at \(d_{\mathrm H}=0.3\), where the MARE is \(4.11\%\) for VAE-MLP and \(7.97\%\) for Descriptor-MLP.
At this level, Descriptor-MLP also shows a stronger positive bias and broader dispersion in signed relative error.
The corresponding aggregate metrics are reported in Table~\ref{tab:e3_prediction_performance}.

\begin{figure}[htbp]
    \centering
    \includegraphics[width=1.0\linewidth]{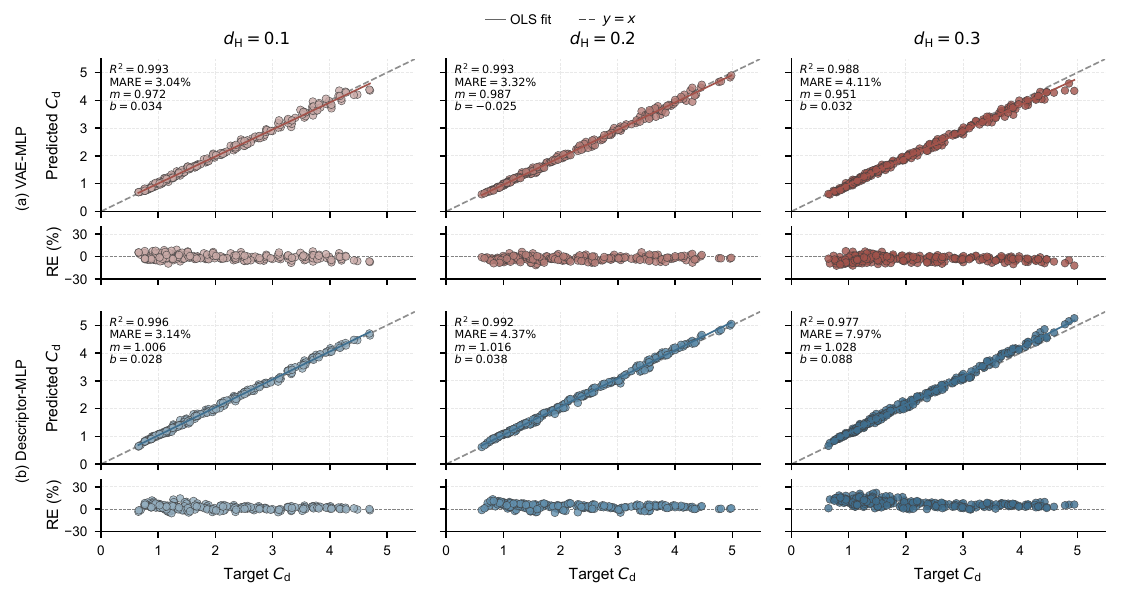}
    \caption{Drag predictions and signed relative errors for training seed 46 at the unseen high-order perturbation levels in E3. Columns correspond to \(d_{\mathrm{H}}=0.1\), \(0.2\), and \(0.3\), while rows correspond to (a) VAE-MLP and (b) Descriptor-MLP. Solid colored lines denote panel-specific ordinary least-squares fits, and gray dashed lines denote ideal parity. Each panel reports \(R^2\), MARE, fitted slope \(m\), and intercept \(b\).}
    \label{fig:e3_parity}
\end{figure}

\begin{table}[htbp]
    \centering
    \small
    \setlength{\tabcolsep}{4pt}
    \caption{Aggregate predictive performance in E3 for training seed 46 under unseen high-order perturbations.}
    \label{tab:e3_prediction_performance}
    \begin{tabular}{lccccc}
        \toprule
        Model
        & \(R^2\)
        & MARE (\%)
        & Max. \(\varepsilon_{\mathrm{rel}}\) (\%)
        & \multicolumn{2}{c}{cMARE (\%)} \\
        \cmidrule(lr){5-6}
        & & & & Minimum & Maximum \\
        \midrule
        VAE-MLP
        & \textbf{0.991}
        & \textbf{3.49}
        & \textbf{12.37}
        & \textbf{1.19}
        & \textbf{6.75} \\
        Descriptor-MLP
        & 0.988
        & 5.16
        & 22.01
        & 1.43
        & 13.80 \\
        \bottomrule
    \end{tabular}
\end{table}

\begin{figure}[htbp]
    \centering
    \includegraphics[width=0.8\linewidth]{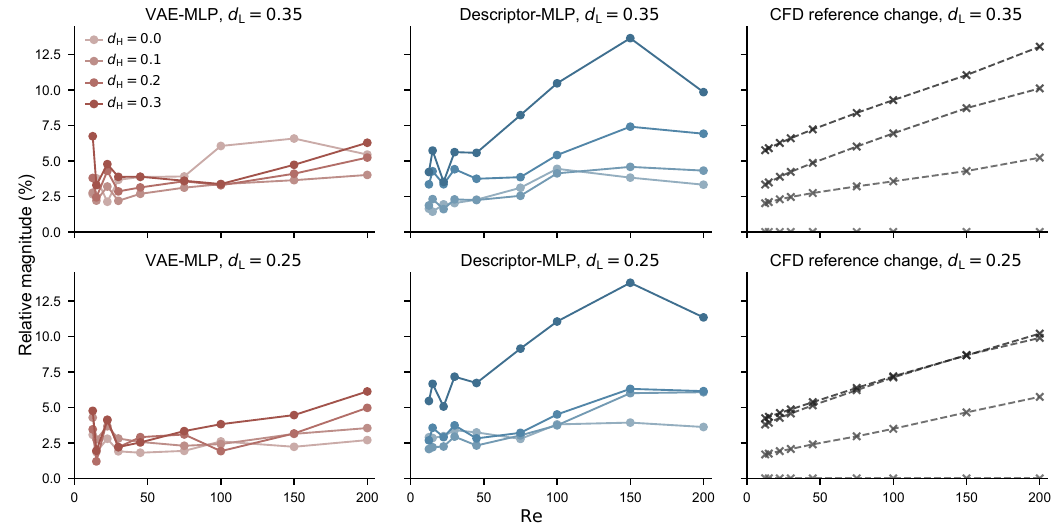}
    \caption{E3 MARE and CFD-reference change as functions of Reynolds number, grouped by \(d_{\mathrm{L}}\). Solid curves denote model MARE, and dashed curves denote the CFD-reference change. Model results are shown for training seed 46.}
    \label{fig:e3_mare_vs_re_by_dl}
\end{figure}

\begin{figure}[htbp]
    \centering
    \includegraphics[width=0.95\linewidth]{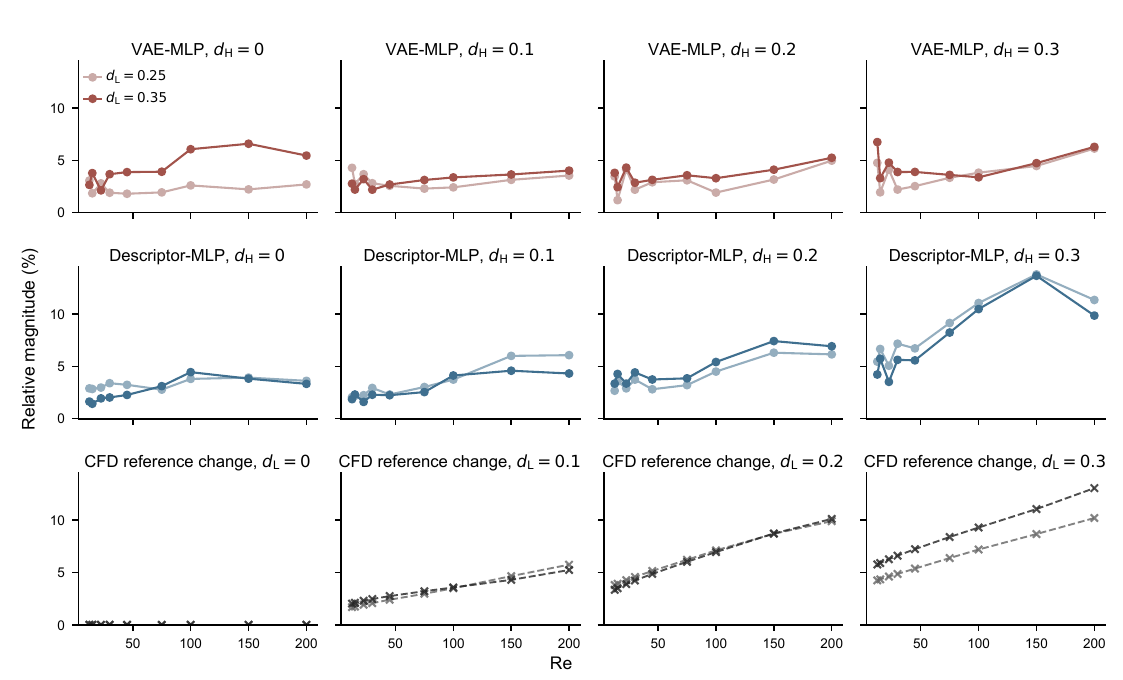}
    \caption{Reynolds-number-resolved E3 results grouped by \(d_{\mathrm{H}}\). Solid curves denote model MARE, and dashed curves denote the CFD-reference change. Model results are shown for training seed 46.}
    \label{fig:e3_mare_vs_re_by_dh}
\end{figure}

For training seed 46, Figures~\ref{fig:e3_mare_vs_re_by_dl} and~\ref{fig:e3_mare_vs_re_by_dh} show model MARE as a function of Reynolds number, grouped by \(d_{\mathrm L}\) and \(d_{\mathrm H}\), respectively.
Over most of the evaluated Reynolds-number range, MARE varies less across \(d_{\mathrm H}\) levels for VAE-MLP than for Descriptor-MLP.
At \(d_{\mathrm H}=0.3\), the MARE gap between the models is pronounced at intermediate-to-high Reynolds numbers and narrows at \(\mathrm{Re}=200\).

\begin{figure}[htbp]
    \centering
    \includegraphics[width=0.95\linewidth]{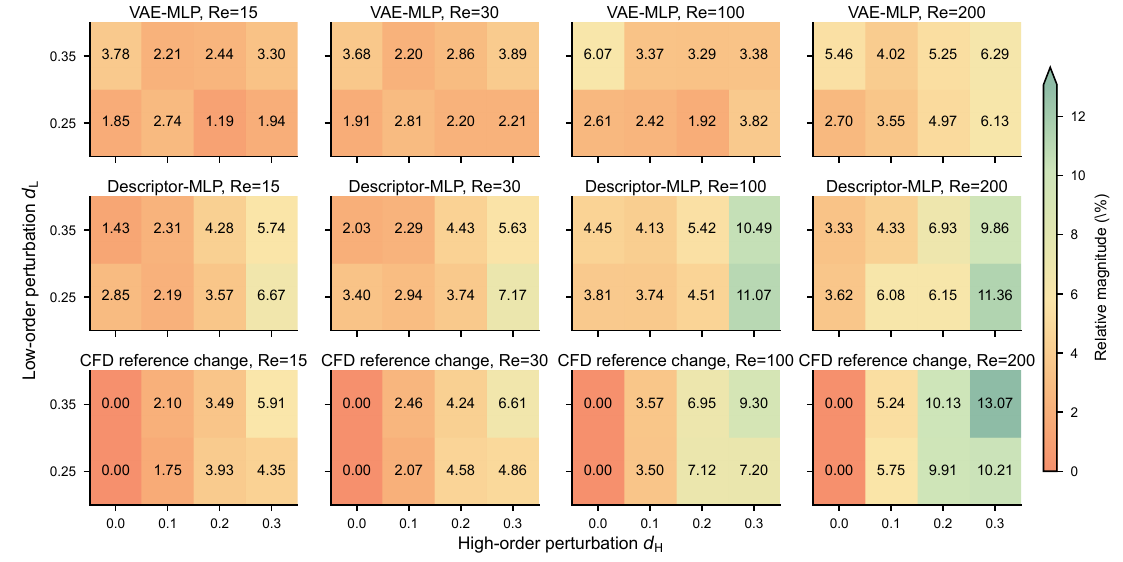}
    \caption{E3 heatmaps over the \((d_{\mathrm L},d_{\mathrm H})\) parameter space at selected Reynolds numbers: VAE-MLP MARE (top), Descriptor-MLP MARE (middle), and CFD-reference change (bottom). Model results are shown for training seed 46.}
    \label{fig:e3_mare_heatmap}
\end{figure}

The CFD-reference change is the mean absolute relative difference between CFD drag coefficients and their parameter-matched \(d_{\mathrm H}=0\) reference values. It generally increases with perturbation level and Reynolds number, as shown by the dashed curves in Figures~\ref{fig:e3_mare_vs_re_by_dl} and~\ref{fig:e3_mare_vs_re_by_dh}. At \(\mathrm{Re}=100\) and \(d_{\mathrm L}=0.25\), however, Figure~\ref{fig:e3_mare_heatmap} shows that the CFD-reference change rises only from \(7.12\%\) to \(7.20\%\) as \(d_{\mathrm H}\) increases from 0.2 to 0.3, while Descriptor-MLP cMARE rises from \(4.51\%\) to \(11.07\%\). At this condition, Descriptor-MLP error therefore increases sharply while the CFD-reference change remains nearly unchanged.

\subsubsection{Analysis of extrapolative performance}
\label{sec:e3_extrapolative_analysis}

To characterize the different extrapolative responses, we compare the predicted and CFD drag differences between geometries at \(d_{\mathrm H}>0\) and matched references at \(d_{\mathrm H}=0\), holding Reynolds number, aspect ratio, incident angle, and \(d_{\mathrm L}\) fixed. We refer to these differences between high-order perturbation levels as cross-level drag differences, defined as
\begin{equation}
    \begin{split}
        \Delta C_{\mathrm d}(d_{\mathrm H})
        & = C_{\mathrm d}(d_{\mathrm H})-C_{\mathrm d}(0),           \\
        \Delta\hat C_{\mathrm d}(d_{\mathrm H})
        & = \hat C_{\mathrm d}(d_{\mathrm H})-\hat C_{\mathrm d}(0),
    \end{split}
    \label{eq:e3_cross_level_drag_difference}
\end{equation}
for \(d_{\mathrm H}\in\{0.1,0.2,0.3\}\). Geometries at different \(d_{\mathrm H}\) levels are independently generated, so these quantities describe cross-level differences between condition-matched realizations.

Agreement with the CFD cross-level differences is quantified using Pearson correlation, sign agreement, and the mean absolute error of \(\Delta\hat C_{\mathrm d}-\Delta C_{\mathrm d}\), termed increment MAE. Sign agreement is the percentage of matched observations for which the predicted and CFD differences have the same sign. The aggregate statistics are calculated over 648 matched observations within each training seed and then summarized across seeds.

Table~\ref{tab:e3_multiseed_response} summarizes the comparison across five training runs. Descriptor-MLP shows higher correlation with the CFD cross-level differences (\(0.733\) versus \(0.577\)) and higher sign agreement (\(90.59\%\) versus \(64.91\%\)), while VAE-MLP has lower increment MAE (\(0.0876\) versus \(0.1025\)). The lower absolute difference error of VAE-MLP therefore accompanies weaker agreement in response direction.

\begin{figure}[htbp]
    \centering
    \includegraphics[width=0.95\linewidth]{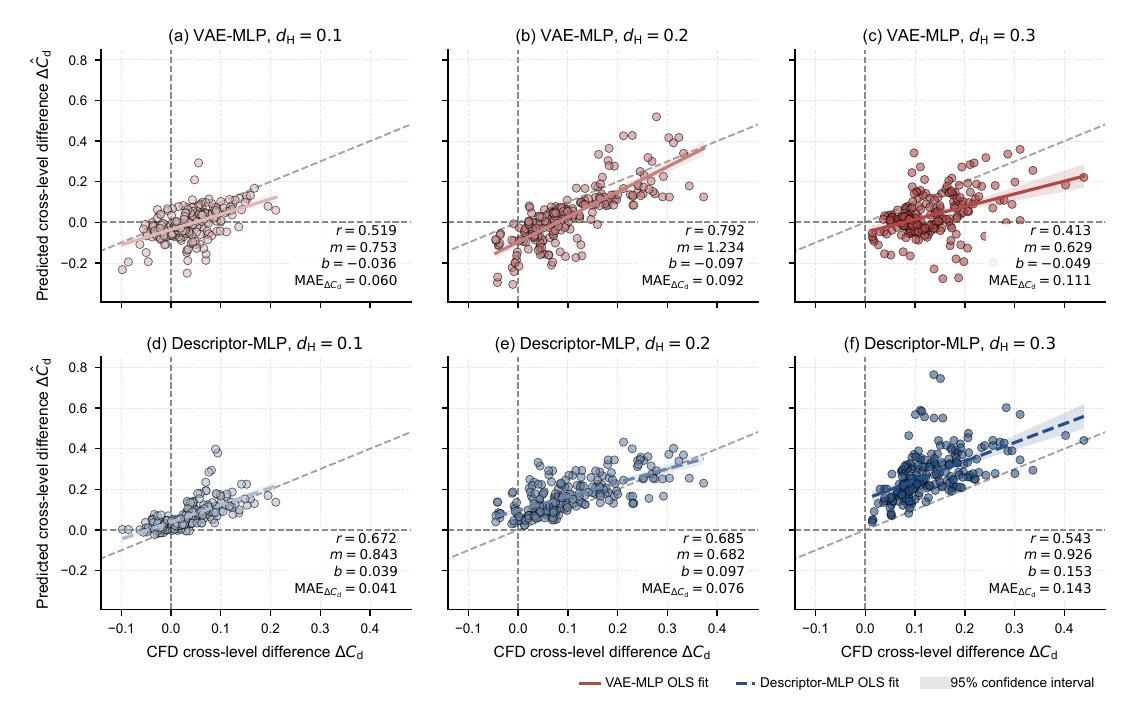}
    \caption{CFD and seed-46 predicted drag differences for VAE-MLP (top) and Descriptor-MLP (bottom) at \(d_{\mathrm H}=0.1\), \(0.2\), and \(0.3\). Each point represents one matched flow observation comparing independently generated geometries at \(d_{\mathrm H}=0\) and the indicated \(d_{\mathrm H}>0\), with aspect ratio, incident angle, \(d_{\mathrm L}\), and Reynolds number held fixed. Dashed lines indicate parity; solid lines and shaded bands show ordinary least-squares fits with 95\% confidence intervals. Annotations report Pearson \(r\), slope, intercept, and increment MAE.}
    \label{fig:e3_cross_level_response}
\end{figure}

At \(d_{\mathrm H}=0.3\), the cross-seed mean predicted drag difference is \(0.3106\) for Descriptor-MLP and \(0.0346\) for VAE-MLP, compared with \(0.1313\) for CFD. Descriptor-MLP therefore overpredicts the mean cross-level change, while VAE-MLP underpredicts it. These opposite deviations distinguish the response magnitudes of the two models despite the lower final-prediction error of VAE-MLP.
Figure~\ref{fig:e3_cross_level_response} illustrates the level-specific responses for seed 46. Descriptor-MLP gives slightly lower increment MAE at \(d_{\mathrm H}=0.1\) and \(0.2\), whereas VAE-MLP gives lower increment MAE at \(d_{\mathrm H}=0.3\). At the strongest perturbation, the mean predicted differences follow the same pattern as the five-seed averages: Descriptor-MLP overpredicts the CFD change, while VAE-MLP underpredicts it.

\begin{table}[htbp]
    \centering
    \small
    \setlength{\tabcolsep}{4pt}
    \caption{Cross-seed summaries of the E3 drag response. The upper block gives mean \(\pm\) sample SD of statistics calculated over all three nonzero levels within each seed. The lower block gives cross-seed means of the level-specific mean predicted increments and the fixed CFD means. Increment MAE and mean increments are in drag-coefficient units.}
    \label{tab:e3_multiseed_response}
    \begin{tabular}{lccc}
        \toprule
        Quantity                              & VAE-MLP             & Descriptor-MLP      & CFD \\
        \midrule
        Pearson \(r\)                         & \(0.577\pm0.098\)   & \(0.733\pm0.037\)   & --- \\
        Sign agreement (\%)                   & \(64.91\pm2.39\)    & \(90.59\pm0.42\)    & --- \\
        Increment MAE                         & \(0.0876\pm0.0059\) & \(0.1025\pm0.0295\) & --- \\
        \midrule
        Mean increment, \(d_{\mathrm H}=0.1\) & \(-0.0062\)         & 0.0738              & 0.0378 \\
        Mean increment, \(d_{\mathrm H}=0.2\) & 0.0378              & 0.1783              & 0.1029 \\
        Mean increment, \(d_{\mathrm H}=0.3\) & 0.0346              & 0.3106              & 0.1313 \\
        \bottomrule
    \end{tabular}
\end{table}

Input-substitution diagnostics assess prediction sensitivity at \(d_{\mathrm H}=0.3\) while keeping each trained predictor and Reynolds number fixed.
All substituted inputs are evaluated against the CFD target of the perturbed geometry. The substitution protocol and full quantitative results are provided in~\ref{app:e3_computational_procedures}. For Descriptor-MLP, restoring \(\Phi\) or \(\Phi_{\mathrm{cross}}\) to its matched reference value reduces the cross-seed mean MARE from \(9.88\%\) to \(3.28\%\) or \(4.50\%\), respectively, whereas restoring \(D_{\mathrm{eq}}\) has little effect. Together with the overpredicted mean cross-level change, these results associate the elevated error with the fitted model's response to changes in the area-related descriptors.

For VAE-MLP, using the latent input of the perturbed geometry gives a cross-seed mean MARE of \(4.97\%\), compared with \(6.34\%\) when the matched reference input is substituted. This improvement occurs in all five training runs. The geometry-specific latent input therefore supplies useful predictive information, even though the model underpredicts the mean cross-level drag change.


\section{Conclusions}
This study develops a VAE-MLP framework that learns a geometry-aware latent representation directly from voxelized particle shapes for drag coefficient prediction. The VAE reconstructs particle geometry from the learned latent representation. The MLP uses the same representation together with the Reynolds number to predict the drag coefficient. The encoder, decoder, and regressor are jointly optimized using reconstruction, KL-regularization, and drag-prediction losses. Comparisons with Descriptor-MLP yield the following findings:

\begin{enumerate}[label=(\arabic*)]
    \item When the evaluated geometric variations are represented during training, VAE-MLP and Descriptor-MLP achieve comparable average prediction accuracy for regular ellipsoids. VAE-MLP achieves lower errors for particles with irregular geometry spanning overall-shape and fine-scale surface variations.

    \item When high-order surface variations are absent from training, VAE-MLP achieves lower overall MARE than Descriptor-MLP in all five training runs with a fixed data split. The overall MARE is \(3.99\%\pm0.55\%\) for VAE-MLP and \(6.11\%\pm1.72\%\) for Descriptor-MLP, expressed as the mean \(\pm\) sample standard deviation across training runs. At the strongest tested perturbation, \(d_{\mathrm H}=0.3\), the corresponding mean MAREs are \(4.97\%\) and \(9.88\%\).

    \item The two models differ in how they use geometric information to predict drag on unseen geometries. At \(d_{\mathrm H}=0.3\), Descriptor-MLP overestimates the mean CFD drag increment, and its elevated prediction error is associated with its fitted response to changes in sphericity and crosswise sphericity. For VAE-MLP, using the latent representations of the perturbed geometries yields lower MARE than using those of the matched references in all five training runs. This improvement demonstrates that the latent inputs provide useful geometric information for predicting drag on the tested unseen geometries.
\end{enumerate}

The binary voxel representation used in this study has limited ability to resolve fine-scale surface features at a fixed spatial resolution. Future work will investigate signed distance functions as an alternative geometric representation to better capture surface details and assess their effects on drag prediction. Due to limited computational resources, the present study focuses exclusively on drag coefficient prediction. Future studies will extend the framework to other hydrodynamic responses, including lift and torque, to further assess its potential as a general geometry-aware prediction approach.

\section*{Acknowledgements}
This work was supported by the JST FOREST Program (Grant Number JPMJFR215S), JSPS KAKENHI (Grant Numbers JP23K26356, JP26K02921 and JP25KJ0026).
Zihao Li acknowledges financial support from the Program of China Scholarship Council through a doctoral scholarship (Grant Number 202406150033).

\section*{CRediT authorship contribution statement}

Zihao Li: Conceptualization, Formal analysis, Data curation, Methodology, Software, Validation, Visualization, Writing – original draft, Writing – review \& editing. Shota Deguchi: Writing – review \& editing, Investigation, Methodology. Mitsuteru Asai: Writing – review \& editing, Funding acquisition, Supervision. Shunhua Chen: Writing – review \& editing, Supervision, Methodology, Resources. Naoto Mitsume: Writing – review \& editing, Conceptualization, Funding acquisition, Supervision, Resources, Project administration, Methodology.

\section*{Declaration of generative AI and AI-assisted technologies in the manuscript preparation process}

During the preparation of this work, the authors used ChatGPT
(OpenAI) to assist with language editing, improving readability,
and refining manuscript organization. After using this tool,
the authors reviewed and edited the content as needed and take
full responsibility for the content of the published article.

\appendix
\section{Spherical-harmonic geometry construction}
\label{app:geometry_generation_details}

The particle surface is represented by separate spherical-harmonic expansions of the three Cartesian coordinates:
\begin{equation}
    \begin{aligned}
        x(\theta,\phi)
        & =
        \sum_{n=0}^{N}\sum_{m=-n}^{n}
        c^{x}_{n,m}Y_n^m(\theta,\phi), \\
        y(\theta,\phi)
        & =
        \sum_{n=0}^{N}\sum_{m=-n}^{n}
        c^{y}_{n,m}Y_n^m(\theta,\phi), \\
        z(\theta,\phi)
        & =
        \sum_{n=0}^{N}\sum_{m=-n}^{n}
        c^{z}_{n,m}Y_n^m(\theta,\phi).
    \end{aligned}
    \label{eq:cartesian_sh_surface}
\end{equation}
Here, \(\theta\) and \(\phi\) are the polar and azimuthal angles, respectively; \(n\) and \(m\) denote the spherical-harmonic degree and order; and \(Y_n^m\) is the corresponding basis function. The coefficients \(c^{x}_{n,m}\), \(c^{y}_{n,m}\), and \(c^{z}_{n,m}\) are complex-valued. The expansion is truncated at \(N=15\), and the real parts of the expanded coordinates are used to obtain the surface vertices. The base ellipsoid provides the coefficients defining the underlying spheroidal shape, while additional coefficients introduce the prescribed geometric irregularity.

The low- and high-order perturbation bands correspond to degrees \(n=2,\ldots,8\) and \(n=9,\ldots,14\), respectively, and are controlled by the parameters \(d_{\mathrm{L}}\) and \(d_{\mathrm{H}}\) defined in Section~\ref{sec:geometry_generation}. Let \(\bm F\) denote the spherical-harmonic coefficient matrix of the volume-normalized base ellipsoid. Its Frobenius norm defines the reference coefficient scale,
\begin{equation}
    d_1
    =
    \|\bm{F}\|_{\mathrm{F}}
    =
    \left(
        \sum_{i,j}
        |F_{ij}|^2
    \right)^{1/2}.
\end{equation}
For the base ellipsoid, the principal-axis scaling factors relative to the unit-diameter reference sphere are
\begin{equation}
    a_{\mathrm{norm}}=\mathrm{Ar}^{2/3},
    \qquad
    b_{\mathrm{norm}}=c_{\mathrm{norm}}=\mathrm{Ar}^{-1/3}.
\end{equation}
Their product is unity, so variation of \(\mathrm{Ar}\) preserves the volume of the unit-diameter reference sphere. For this ellipsoid, the corresponding reference scale is
\begin{equation}
    d_1
    =
    \sqrt{
        \frac{\pi}{3}
        \left(
            \mathrm{Ar}^{4/3}
            +
            2\mathrm{Ar}^{-2/3}
        \right)
    }.
\end{equation}
Thus, \(d_1\) varies with the aspect ratio of the base ellipsoid and provides the reference scale for the prescribed perturbations.

The parameters \(d_{\mathrm L}\) and \(d_{\mathrm H}\) are converted to intermediate coefficient scales according to
\begin{equation}
    D_{\mathrm{L}}
    =
    d_{\mathrm{L}}d_1
    \left[
        \sum_{n=2}^{8}
        \left(\frac{2}{n}\right)^{1.387}
    \right]^{-1},
    \qquad
    D_{\mathrm{H}}
    =
    d_{\mathrm{H}}d_1
    \left[
        \sum_{n=9}^{15}
        \left(\frac{9}{n}\right)^{1.426}
    \right]^{-1}.
\end{equation}
In the implemented geometry generator, the coefficient norm assigned to each Cartesian component at degree \(n\) is
\begin{equation}
    q_n=
    \begin{cases}
        D_{\mathrm L},
        & n=2,           \\[5pt]
        \dfrac{D_{\mathrm L}}{\sqrt{3}}
        \left(\dfrac{2}{n-1}\right)^{1.387},
        & n=3,\ldots,8,  \\[9pt]
        \dfrac{D_{\mathrm H}}{\sqrt{3}}
        \left(\dfrac{9}{n-1}\right)^{1.426},
        & n=9,\ldots,14, \\[9pt]
        0,
        & n=15.
    \end{cases}
\end{equation}

For each degree \(n\) and Cartesian coordinate \(j\in\{x,y,z\}\), let \(\bm p_n^j\) denote the randomly generated complex coefficient vector over the SH orders used by the implementation. Each Cartesian coefficient block is normalized independently and rescaled as
\begin{equation}
    \bm c_n^j
    =
    q_n
    \frac{\bm p_n^j}
    {\|\bm p_n^j\|_2},
    \qquad
    j\in\{x,y,z\}.
\end{equation}
Consequently, the three Cartesian coefficient blocks have the same coefficient norm \(q_n\) at a given degree.

A new random coefficient realization is generated independently for each particle shape. The resulting coefficients are used in the Cartesian SH expansions of Eq.~\eqref{eq:cartesian_sh_surface}, and the real parts of the expanded coordinates are retained to construct the surface vertices. Different \(d_{\mathrm L}\) and \(d_{\mathrm H}\) levels therefore correspond to independently generated particle realizations rather than to nested perturbations obtained by rescaling a common coefficient realization.

For surface discretization, a unit-sphere mesh is constructed from a regular icosahedron. Its triangular faces are subdivided, and the resulting vertices are projected onto the unit sphere. The vertex directions are then converted to \((\theta,\phi)\) and evaluated using Eq.~\eqref{eq:cartesian_sh_surface}. The triangular connectivity of the subdivided sphere is retained after the SH mapping, resulting in a surface mesh with 5120 triangular facets for the mesh resolution used in this study. The SH perturbations are applied to the nominally volume-normalized base ellipsoid without rescaling the completed intrinsic geometry to a common volume. The prescribed incident angle is imposed after perturbation. The implementation of this geometry-generation procedure is based on the core particle-generation routines of SHPSG~\cite{zhao2020shpsg}.

\section{CFD formulation and reference-data generation}

This section provides additional details of the CFD formulation and numerical procedure summarized in Section~\ref{sec:cfd_reference_data}. The simulations are performed in OpenFOAM v2306~\cite{weller1998tensorial} for a rigid particle held at a prescribed position and orientation in an incompressible Newtonian flow.

\subsection{Governing equations and computational setup}
\label{app:cfd_theory}

The flow field is governed by the incompressible Navier--Stokes equations,
\begin{equation}
    \nabla \cdot \bm{U} = 0,
    \label{eq:appendix_continuity}
\end{equation}
and
\begin{equation}
    \frac{\partial \bm{U}}{\partial t}
    +
    \nabla \cdot
    \left(
        \bm{U}\otimes\bm{U}
    \right)
    =
    -\nabla p_k
    +
    \nu \nabla^2 \bm{U},
    \label{eq:appendix_momentum}
\end{equation}
where \(\bm{U}\) is the fluid velocity, \(\nu\) is the kinematic viscosity, and \(p_k=p/\rho_\infty\) denotes the kinematic pressure used in the incompressible formulation.

The computational domain extends from \(-10D_{\mathrm{ref}}\) to \(20D_{\mathrm{ref}}\) in the streamwise direction and from \(-10D_{\mathrm{ref}}\) to \(10D_{\mathrm{ref}}\) in both transverse directions. A uniform velocity
\(\bm{U}=(1,0,0)\) is prescribed at the inlet. The outlet uses a fixed-pressure condition together with a pressure-inlet--outlet velocity condition. Slip conditions are applied at the lateral boundaries, while the particle surface is treated as a stationary no-slip wall.

The free-stream velocity is fixed at \(U_\infty=1.0\). The fixed length \(D_{\mathrm{ref}}=1.0\) is used only to specify the computational-domain dimensions above. For each final particle geometry, the volume-equivalent diameter \(D_{\mathrm{eq}}=(6V/\pi)^{1/3}\) is calculated from its volume, and the prescribed Reynolds number is controlled through the kinematic viscosity according to
\begin{equation}
    \nu
    =
    \frac{U_\infty D_{\mathrm{eq}}}{\mathrm{Re}}.
    \label{eq:appendix_reynolds_control}
\end{equation}
The simulations cover \(\mathrm{Re}\in[10,200]\). A common reference density \(\rho_\infty=1\) is used. Each particle uses the geometry-specific drag-reference area defined in Eq.~\eqref{eq:reference_area}.

\subsection{Numerical solution and drag evaluation}
\label{app:cfd_drag_evaluation}
The transient incompressible equations are solved using the OpenFOAM v2306 \texttt{pimpleFoam} solver. A first-order Euler scheme is used for time integration. Spatial discretization uses a Gauss linear scheme for gradients, Gauss linearUpwind with the velocity gradient for the velocity-convection term, and a Gauss linear corrected scheme for Laplacian terms. Pressure and velocity are coupled using the PIMPLE algorithm~\cite{greenshields2022cfd}. Five pressure--velocity correctors and two non-orthogonal correctors are applied within each time step. Adaptive time stepping is used with the maximum Courant number limited to 1.

The hydrodynamic force on the particle is evaluated from the pressure and skin-friction contributions integrated over the particle-surface patches. The resulting force is projected onto the prescribed drag direction, and the corresponding drag coefficient is obtained using the reference free-stream velocity, density, and area~\cite{openfoam_forceCoeffs}. This procedure is consistent with the definition of \(C_{\mathrm{d}}\) given in Eq.~\eqref{eq:drag_coefficient_definition}.

For each simulation, the drag coefficient is recorded at every computational time step, and the scalar learning target is obtained by arithmetic averaging over a terminal subset of the recorded \(C_{\mathrm{d}}\) samples. The averaging window comprises the final 20\% of the recorded samples, with a minimum of 50 samples whenever available. If \(N\) drag coefficient samples are available, the number of samples included in the averaging window is
\begin{equation}
    N_{\mathrm{avg}}
    =
    \min
    \left[
        N,\,
        \max
        \left(
            \lfloor 0.2N \rfloor,
            50
        \right)
    \right],
    \label{eq:appendix_cd_window}
\end{equation}
and the reference value is
\begin{equation}
    \overline{C}_{\mathrm{d}}
    =
    \frac{1}{N_{\mathrm{avg}}}
    \sum_{i=N-N_{\mathrm{avg}}+1}^{N}
    C_{\mathrm{d},i}.
    \label{eq:appendix_cd_average}
\end{equation}
Across the assembled CFD datasets, the relative standard deviation of the selected terminal samples has median and 95th-percentile values of \(0.0035\%\) and \(0.152\%\), respectively. The corresponding 95th percentile remains below \(0.3\%\) within each dataset, indicating low terminal variability for most simulations.

\subsection{CFD mesh verification}
\label{app:cfd_mesh_verification}

A three-level mesh study assesses the production-grid resolution. The coarse, medium, and fine meshes contain approximately \(43\times10^3\), \(154\times10^3\), and \(424\times10^3\) cells, respectively. Local refinement is applied near the particle surface and in the downstream wake, together with six surface layers using an expansion ratio of 1.2.

\begin{figure}[htbp]
    \centering
    \includegraphics[width=0.6\linewidth]{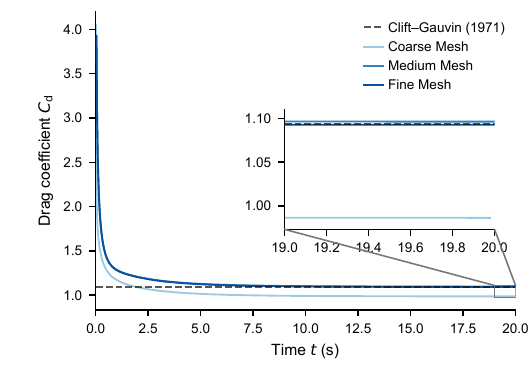}
    \caption{Time histories of the drag coefficient for a spherical particle obtained using the coarse, medium, and fine meshes. The dashed line denotes the sphere correlation of Clift--Gauvin~\cite{clift1971motion}.}
    \label{fig:cd_mesh_convergence}
\end{figure}

Figure~\ref{fig:cd_mesh_convergence} compares the resulting drag coefficient histories for a spherical particle. The medium- and fine-grid solutions approach similar values over the final portion of the simulations, whereas a larger difference is observed on the coarse grid.

The mesh assessment also includes representative oblate, prolate, and spherical particles spanning different aspect ratios, incident angles, and Reynolds numbers. The relative difference between the medium- and fine-grid predictions is defined as
\begin{equation}
    \delta_{\mathrm{M/F}}
    =
    \frac{
        \left|
        C_{\mathrm{d}}^{\mathrm{M}}
        -
        C_{\mathrm{d}}^{\mathrm{F}}
        \right|
    }{
        C_{\mathrm{d}}^{\mathrm{F}}
    }
    \times 100\%.
    \label{eq:appendix_mesh_difference}
\end{equation}

\begin{table}[htbp]
    \centering
    \caption{Medium-to-fine mesh differences for representative particle conditions.}

    \begin{tabular}{ccccccc}
        \toprule
        Case                                   &
        \(\mathrm{Ar}\)                        &
        \(\alpha_{\mathrm{inc}}\) (\(^\circ\)) &
        \(\mathrm{Re}\)                        &
        \(C_{\mathrm{d}}^{\mathrm{M}}\)        &
        \(C_{\mathrm{d}}^{\mathrm{F}}\)        &
        \(\delta_{\mathrm{M/F}}\) (\%)
        \\
        \midrule
        1                                      & 0.25 & 90 & 200 & 0.72354 & 0.72237 & 0.162 \\
        2                                      & 0.25 & 45 & 100 & 1.93353 & 1.92132 & 0.636 \\
        3                                      & 0.25 & 0  & 10  & 7.32888 & 7.30070 & 0.386 \\
        4                                      & 2.5  & 90 & 200 & 1.15412 & 1.14819 & 0.516 \\
        5                                      & 2.5  & 45 & 100 & 1.17939 & 1.17612 & 0.278 \\
        6                                      & 2.5  & 0  & 10  & 3.85681 & 3.86200 & 0.134 \\
        7                                      & 1.0  & -- & 200 & 0.78027 & 0.77609 & 0.539 \\
        8                                      & 1.0  & -- & 100 & 1.09574 & 1.09274 & 0.275 \\
        \bottomrule
    \end{tabular}%

\end{table}

Across these representative conditions, the medium-to-fine differences remain below \(0.64\%\), with a mean difference of \(0.366\%\). Together with the convergence behavior shown in Figure~\ref{fig:cd_mesh_convergence}, these results support the numerical framework and the use of the medium-mesh strategy for the regular reference particles examined here while substantially reducing the computational cost relative to the fine grid. Because this verification does not include high-order perturbations, it is not a direct demonstration of mesh independence for every irregular geometry in the dataset.

\section{Detailed dataset and evaluation configurations}
\label{app:dataset_configuration}

This appendix specifies the training parameter domains, sample counts, and train--validation configurations for Experiments E1--E3, together with the separate geometry-only dataset used for E1 latent-space visualization. The external evaluation grids are reported in Table~\ref{tab:external_config} in the main text.

\paragraph{Training dataset}
Training geometries are generated using Latin hypercube sampling (LHS) within the particle parameters specified in Table~\ref{tab:training_config}. For each intrinsic geometry, the incident angle and Reynolds number are jointly sampled using a two-dimensional LHS over \(\alpha_{\mathrm{inc}}\in[0^\circ,90^\circ]\) and \(\mathrm{Re}\in[10,200]\). Dataset partitioning is performed at the intrinsic-geometry level, such that all orientations and flow conditions derived from the same intrinsic geometry are assigned to either the training or validation subset.

\begin{table}[htbp]
    \centering
    \caption{Training dataset configurations for Experiments E1--E3. Flow conditions are sampled using a two-dimensional LHS over \(\alpha_{\mathrm{inc}}\in[0^\circ,90^\circ]\) and \(\mathrm{Re}\in[10,200]\).}
    \label{tab:training_config}
    \small
    \begin{tabular}{lcccccc}
        \toprule
        \multirow{2}{*}{\textbf{Exp.}}
            & \multicolumn{3}{c}{\textbf{Particle Parameter Domain}}
            & \multirow{2}{*}{\textbf{Geometries}}
            & \textbf{Flow Conditions}
            & \textbf{Train / Val.} \\
        \cmidrule(lr){2-4}
            & \(\mathrm{Ar}\)
            & \(d_{\mathrm{L}}\)
            & \(d_{\mathrm{H}}\)
            &
            & \textbf{per Geometry}
            & \textbf{Samples} \\
        \midrule
        E1  & \([0.25,2.5]\) & 0             & 0           & 150 & 8 & 1024 / 176 \\
        E2a & \([0.25,2.5]\) & \([0.1,0.4]\) & 0           & 300 & 4 & 1020 / 180 \\
        E2b & \([0.25,2.5]\) & \([0.1,0.4]\) & \([0,0.3]\) & 300 & 4 & 1020 / 180 \\
        E3  & \([0.25,2.5]\) & \([0.1,0.4]\) & 0           & 300 & 4 & 1020 / 180 \\
        \bottomrule
    \end{tabular}
\end{table}

\paragraph{E1 latent-space visualization dataset}
\label{app:e1_latent_visualization_dataset}

For the latent-space analysis in E1, a separate geometry-only dataset is constructed to systematically probe variations in the learned representation with respect to aspect ratio and incident angle. The dataset is defined as the Cartesian product of 15 aspect-ratio levels logarithmically spaced over \(\mathrm{Ar}\in[0.5,2.0]\) and 10 incident-angle levels uniformly spaced over \(\alpha_{\mathrm{inc}}\in[0^\circ,90^\circ]\), yielding 150 geometries in total. These geometries are voxelized using the same preprocessing procedure as the E1 training data. The dataset is used exclusively for latent-space visualization and contains no CFD-derived drag targets.

\section{Neural-network implementation details}

\subsection{Architecture specification}
\label{app:architecture_specification}

\begin{figure}[htbp]
    \centering
    \includegraphics[width=1.0\linewidth]{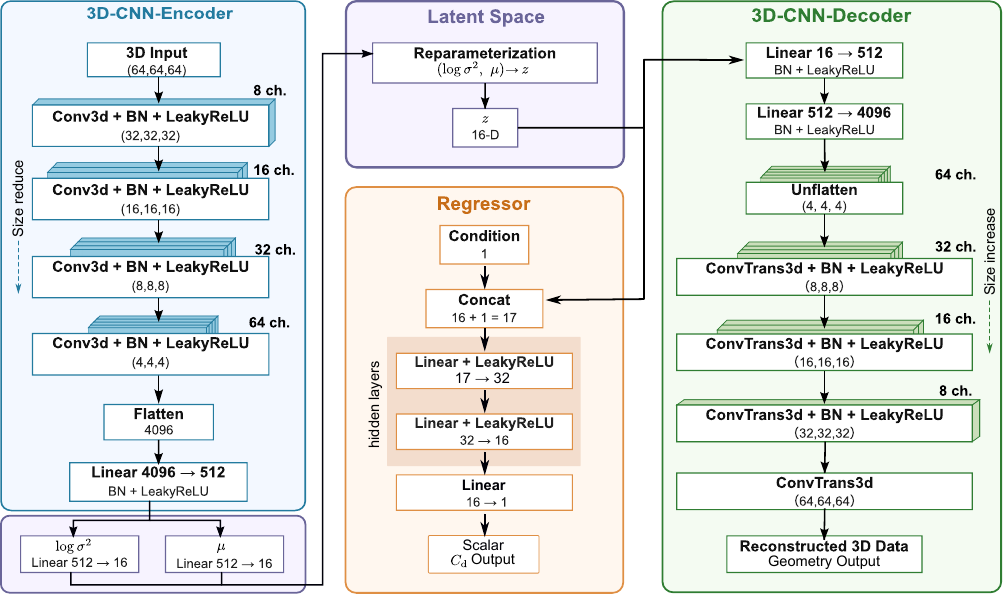}
    \caption{VAE-MLP architecture and tensor-dimensional flow for the experiments using \(d=16\).}
    \label{fig:model_architecture}
\end{figure}

The layer-level implementation of the VAE-MLP is documented in Figure~\ref{fig:model_architecture} and Table~\ref{tab:network_architecture}. The figure shows the complete tensor-dimensional flow, whereas the table reports the operations, channel or unit widths, convolution parameters, output dimensions, normalization, and activation functions. The same layer configuration is used in all main experiments, with \(d=2\) for E1 and \(d=16\) for E2a, E2b, and E3. For layers using batch normalization, BN is applied before the LeakyReLU activation. All LeakyReLU activations use a negative slope of 0.01. The hidden layers of the drag regressor follow the sequence Linear--LeakyReLU. The regressor receives the latent representation concatenated with the Reynolds number, giving an input dimension of \(d+1\), while the final decoder layer produces voxel logits without batch normalization or an output activation.

\begin{table}[htbp]
    \centering
    \small
    \caption{Layer configurations of the VAE-MLP architecture. Conv3D and ConvT3D denote three-dimensional convolution and transposed convolution. In the \((k/s/p)\) column, \(k\), \(s\), and \(p\) denote kernel size, stride, and padding, respectively; RP denotes the reparameterization operation.}
    \label{tab:network_architecture}
    \begin{tabular}{llccccc}
        \toprule
        Module
        & Operation
        & Channels / units
        & \((k/s/p)\)
        & Output
        & BN
        & Activation \\
        \midrule

        \multirow{3}{*}{Encoder}
        & \(4\times\) Conv3D
        & \(8,\,16,\,32,\,64\)
        & \(3/2/1\)
        & \(64\times4^3\)
        & \checkmark
        & LeakyReLU \\

        & Flatten
        & --
        & --
        & \(4096\)
        & --
        & -- \\

        & \(1\times\) Linear
        & \(512\)
        & --
        & \(512\)
        & \checkmark
        & LeakyReLU \\
        \midrule

        \multirow{3}{*}{Latent}
        & Linear (\(\bm{\mu}\))
        & \(d\)
        & --
        & \(\bm{\mu}\in\mathbb{R}^{d}\)
        & --
        & -- \\

        & Linear (\(\bm{v}\))
        & \(d\)
        & --
        & \(\bm{v}\in\mathbb{R}^{d}\)
        & --
        & -- \\

        & RP
        & --
        & --
        & \(\bm{z}\in\mathbb{R}^{d}\)
        & --
        & -- \\
        \midrule

        \multirow{4}{*}{Decoder}
        & \(2\times\) Linear
        & \(512,\,4096\)
        & --
        & \(4096\)
        & \checkmark
        & LeakyReLU \\

        & Unflatten
        & --
        & --
        & \(64\times4^3\)
        & --
        & -- \\

        & \(3\times\) ConvT3D
        & \(32,\,16,\,8\)
        & \(4/2/1\)
        & \(8\times32^3\)
        & \checkmark
        & LeakyReLU \\

        & \(1\times\) ConvT3D
        & \(1\)
        & \(4/2/1\)
        & \(1\times64^3\)
        & --
        & -- \\
        \midrule

        \multirow{2}{*}{Regressor}
        & \(2\times\) Linear
        & \(32,\,16\)
        & --
        & \(16\)
        & --
        & LeakyReLU \\

        & \(1\times\) Linear
        & \(1\)
        & --
        & \(\hat{C}_{\mathrm{d}}\)
        & --
        & -- \\

        \bottomrule
    \end{tabular}
\end{table}

\subsection{Empirical drag correlations}
\label{app:empirical_baseline_formulations}

The ellipsoid-specific correlation of Ke et al.~\cite{ke2018} is evaluated as
\begin{align}
    C_{\mathrm{d}}^{\mathrm{Ke}}
    ={} &
    \frac{c_1}{\mathrm{Re}}\mathrm{Ar}^{c_2}
    +
    \frac{c_3}{\sqrt{\mathrm{Re}}}\mathrm{Ar}^{c_4}
    +
    c_5\mathrm{Ar}^{c_6}
    \nonumber \\
    & +
    \mathrm{Ar}^{c_7}(\mathrm{Ar}-1)
    \frac{c_8}{\mathrm{Re}^{c_9}}
    \sin^2(c_{10}\zeta),
\end{align}
where \(\zeta=(\pi/180)\alpha_{\mathrm{inc}}\) converts the prescribed incident angle \(\alpha_{\mathrm{inc}}\) from degrees to radians. The coefficients are
\begin{equation}
    \begin{aligned}
        (c_1,c_2,c_3,c_4,c_5)
        & =
        (18.7371,\,0.2883,\,7.9738,\,-0.5126,\,0.1938), \\
        (c_6,c_7,c_8,c_9,c_{10})
        & =
        (-1.1848,\,-0.5531,\,2.6334,\,0.2199,\,0.9865).
    \end{aligned}
\end{equation}

The H\"{o}lzer--Sommerfeld correlation~\cite{holzer2008} is evaluated in its crosswise-sphericity form:
\begin{align}
    C_{\mathrm{d}}^{\mathrm{HS}}
    ={} &
    \frac{8}{\mathrm{Re}\sqrt{\Phi_{\mathrm{cross}}}}
    +
    \frac{16}{\mathrm{Re}\sqrt{\Phi}}
    +
    \frac{3}{\sqrt{\mathrm{Re}}\,\Phi^{3/4}}
    \nonumber \\
    & +
    \frac{
        0.42\,
        10^{0.4[-\log_{10}(\Phi)]^{0.2}}
    }{
        \Phi_{\mathrm{cross}}
    }.
\end{align}
Here, \(\Phi\) and \(\Phi_{\mathrm{cross}}\) follow the definitions in Eq.~\eqref{eq:descriptor_definitions}.

\subsection{implementation and random seed detail}
\label{app:training_reproducibility}

The implementation uses PyTorch 2.5.1 with CUDA 12.1 and cuDNN 9.1.0.
A fixed random seed of 42 is used for dataset partitioning in all experiments.
E1, E2a, and E2b are trained using seed 42.
For E3, five independent training runs are performed using seeds 42--46 while retaining the same dataset partition. Results from seed 46 are shown in the detailed E3 figures, whereas the E3 input-substitution tables summarize all five training runs.

\section{Drag variation across geometric perturbations}
\label{app:drag_variation}

\begin{figure}[htbp]
    \centering
    \includegraphics[width=0.6\linewidth]
    {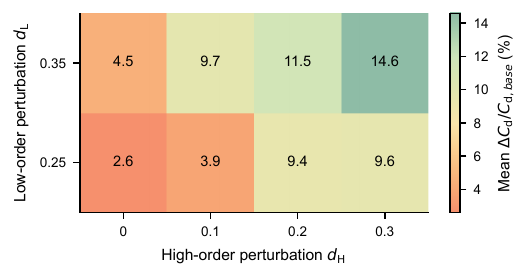}
    \caption{Mean signed relative drag change (\%) from
        parameter-matched regular ellipsoids for different combinations
    of \(d_{\mathrm L}\) and \(d_{\mathrm H}\).}
    \label{fig:drag_variation_heatmap}
\end{figure}

CFD drag coefficients for regular ellipsoids, particles with large-scale shape variations, and particles combining large-scale shape variations and fine-scale surface variations are compared at common parameter combinations: \(\mathrm{Ar}\in\{0.5,0.75,1.5,2.0\}\), \(\alpha_{\mathrm{inc}}\in\{0^\circ,30^\circ,60^\circ\}\), and \(\mathrm{Re}\in\{12.5,15,30,45,100,150,200\}\). Each \((d_{\mathrm L},d_{\mathrm H})\) condition contains 84 samples matched to regular ellipsoids at the same aspect ratio, incident angle, and Reynolds number. For each matched pair, the signed relative drag change is \((C_{\mathrm d}-C_{\mathrm d,ref})/C_{\mathrm d,ref}\), where \(C_{\mathrm d,ref}\) is the CFD target for the matched regular ellipsoid at \(d_{\mathrm L}=d_{\mathrm H}=0\). Figure~\ref{fig:drag_variation_heatmap} shows the mean signed relative drag change, averaged over the 84 matched samples for each perturbation condition and expressed as a percentage. Perturbed geometries at different parameter levels are independently generated realizations.

The sampled irregular-particle conditions yield positive mean drag changes relative to the regular ellipsoids. At \(d_{\mathrm H}=0.3\), the mean increases are \(9.14\%\) and \(13.88\%\) for \(d_{\mathrm L}=0.25\) and \(0.35\), respectively, compared with \(2.57\%\) and \(5.03\%\) at \(d_{\mathrm H}=0\). These results quantify the mean drag changes relative to regular ellipsoids for particles with large-scale shape variations alone and for particles combining large-scale shape variations with fine-scale surface variations.

\section{E2b experiment result supplement}
\label{app:e2b_additional_condition_results}

\begin{figure}[htbp]
    \centering
    \includegraphics[width=1.0\linewidth]{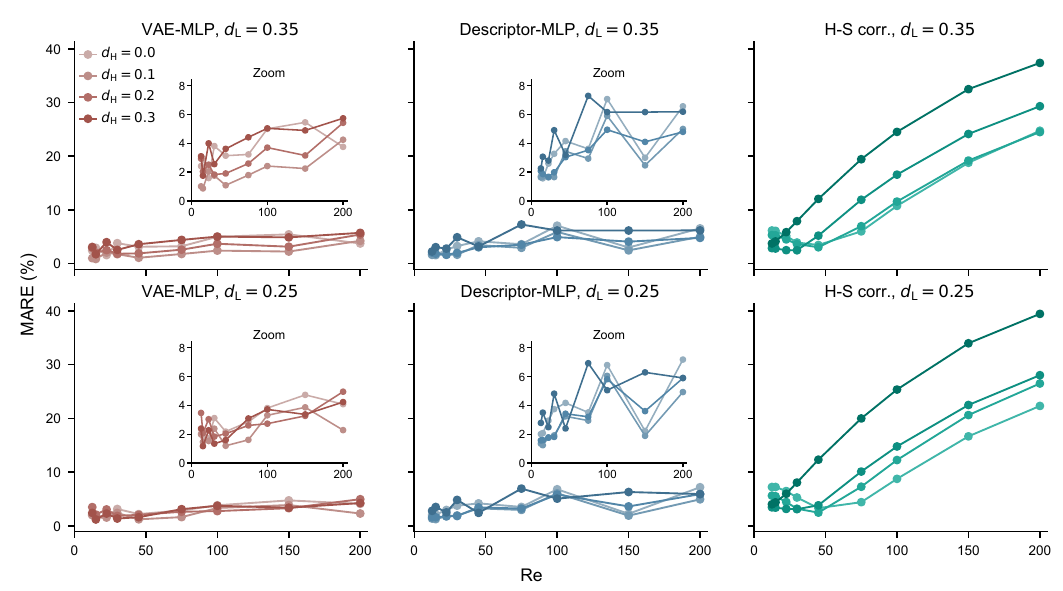}
    \caption{Reynolds-number-resolved MARE in E2b, grouped by \(d_{\mathrm{L}}\).}
    \label{fig:e2b_mare_vs_re_by_dl}
\end{figure}

\begin{figure}[htbp]
    \centering
    \includegraphics[width=1.0\linewidth]{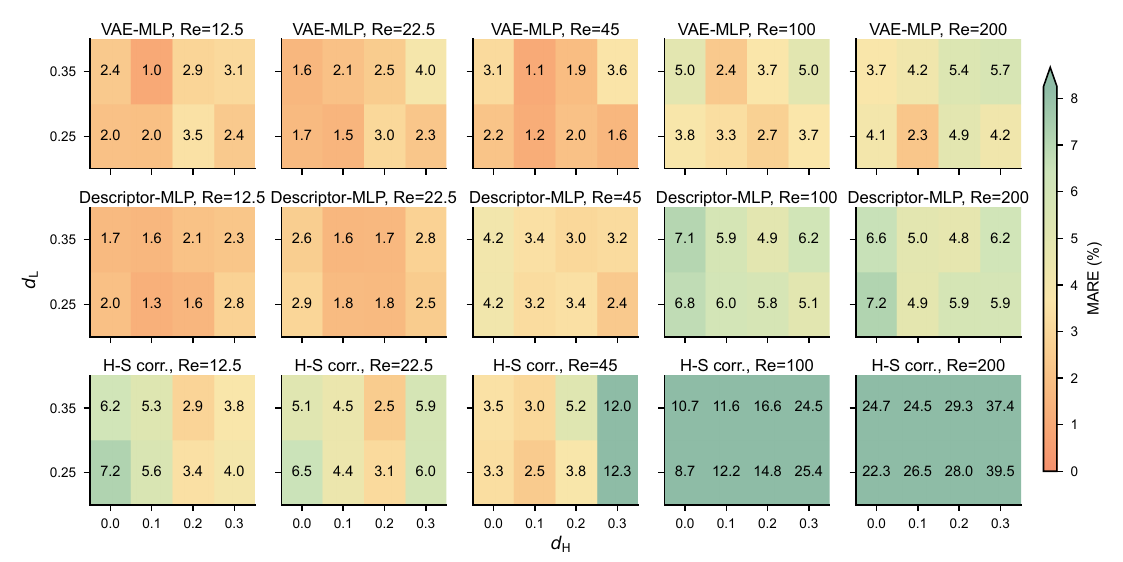}
    \caption{E2b MARE over the \((d_{\mathrm{L}},d_{\mathrm{H}})\) parameter space at the selected Reynolds numbers.}
    \label{fig:e2b_mare_heatmap}
\end{figure}

This appendix provides supplementary results for the E2b experiment, complementing the \(d_{\mathrm{H}}\)-resolved analysis presented in the main text. The results presented in this appendix are based on the same E2b evaluation dataset as that used in the main-text analysis. Figure~\ref{fig:e2b_mare_vs_re_by_dl} provides the complementary grouping by the underlying \(d_{\mathrm{L}}\) level, while Figure~\ref{fig:e2b_mare_heatmap} localizes the errors over the full \((d_{\mathrm{L}},d_{\mathrm{H}})\) evaluation grid at selected Reynolds numbers. The maximum VAE-MLP cMARE of \(5.73\%\) occurs at \(\mathrm{Re}=200\), \(d_{\mathrm{L}}=0.35\), and \(d_{\mathrm{H}}=0.3\). The corresponding maximum for Descriptor-MLP is \(7.28\%\) at \(\mathrm{Re}=75\), \(d_{\mathrm{L}}=0.35\), and \(d_{\mathrm{H}}=0.3\).

\section{Latent sensitivity analysis}
\label{app:latent_sensitivity}
\begin{figure}[htbp]
    \centering
    \includegraphics[width=0.8\linewidth]{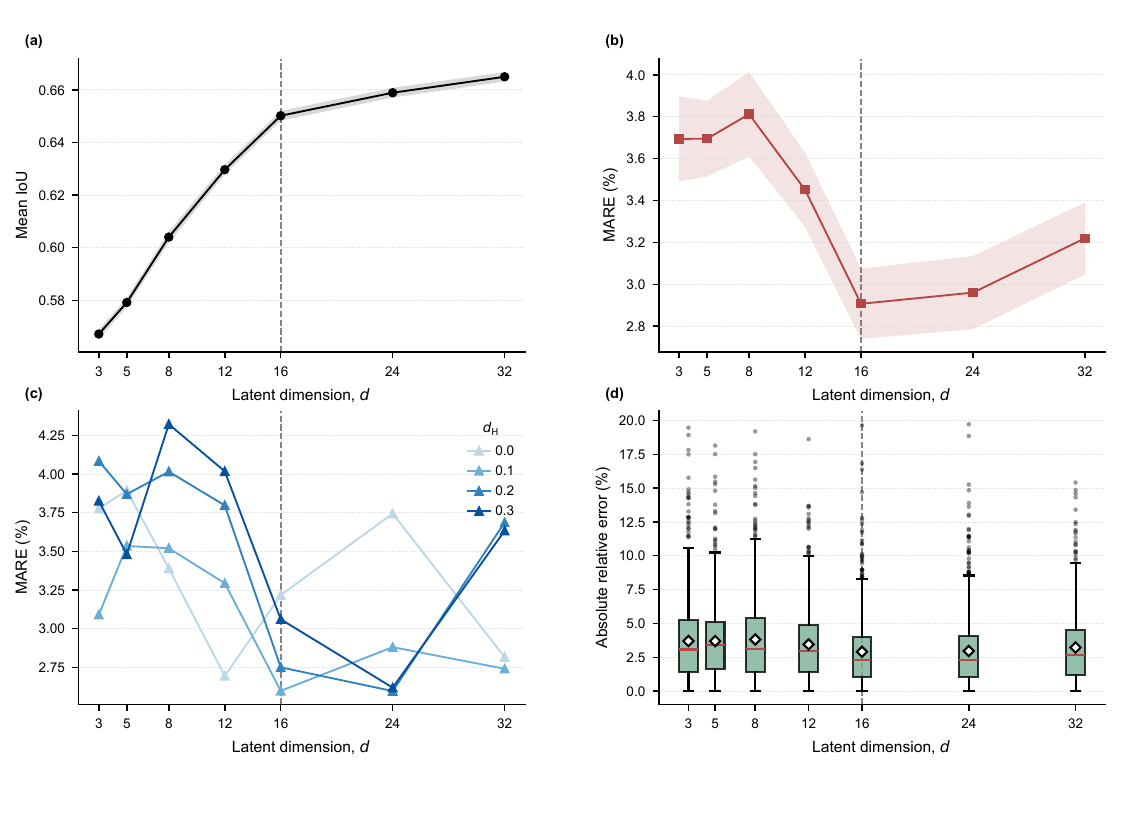}
    \caption{Effects of latent dimension on reconstruction IoU and drag-prediction error in E2b.}
    \label{fig:e2b_latent_sensitivity}
\end{figure}

The latent dimension is varied from 3 to 32 in E2b to assess its effect on reconstruction and drag prediction. Reconstruction is evaluated from the decoder logits. For sample \(i\) and voxel \(v\), the binary reconstructed occupancy \(\hat{x}_{iv}\) and sample-wise intersection-over-union (IoU) are defined as
\begin{equation}
    \hat{x}_{iv}
    =
    \mathbb{I}\!\left[\sigma(\ell_{iv})\geq 0.5\right],
    \qquad
    \operatorname{IoU}_i
    =
    \frac{\sum_v \hat{x}_{iv}x_{iv}}
    {\sum_v\left(\hat{x}_{iv}+x_{iv}-\hat{x}_{iv}x_{iv}\right)},
    \label{eq:reconstruction_iou}
\end{equation}
where \(\ell_{iv}\) is the decoder logit, \(\sigma(\cdot)\) is the sigmoid function, \(x_{iv}\) is the binary target occupancy, and \(\mathbb{I}[\cdot]\) is the indicator function. A value of one is assigned if the union is empty, although this case does not occur for the non-empty particle targets evaluated here. The reported mean IoU is the arithmetic mean of \(\operatorname{IoU}_i\) over the 864 external evaluation samples for each latent dimension. Figure~\ref{fig:e2b_latent_sensitivity}(a) shows that mean IoU increases from approximately 0.57 at \(d=3\) to 0.66 at \(d=32\). The improvement becomes smaller beyond \(d=16\), although reconstruction continues to improve over the evaluated range.

Drag-prediction error follows a different trend.
Overall MARE reaches its minimum at \(d=16\) and increases at \(d=24\) and \(32\), as shown in Table~\ref{tab:e2b_latent_sensitivity}.
The dimension giving the lowest level-specific MARE varies with \(d_{\mathrm H}\): it is \(d=12\) at \(d_{\mathrm H}=0\), \(d=16\) at \(d_{\mathrm H}=0.1\), and \(d=24\) at \(d_{\mathrm H}=0.2\) and \(0.3\).
The sample-wise distributions in Figure~\ref{fig:e2b_latent_sensitivity}(d) likewise narrow from \(d=3\) to 16, with no further improvement in the upper error tail at larger dimensions.

\begin{table}[htbp]
    \centering
    \small
    \caption{Latent-dimension sensitivity results in E2b. IoU is the mean of the sample-wise values defined in Eq.~\eqref{eq:reconstruction_iou} over 864 external evaluation samples.}
    \label{tab:e2b_latent_sensitivity}
    \begin{tabular}{ccccccc}
        \toprule
        \(d\) & Mean IoU & Overall MARE (\%) & \multicolumn{4}{c}{MARE (\%) by \(d_{\mathrm{H}}\)} \\
        \cmidrule(lr){4-7}
              &          &                   & 0.0  & 0.1  & 0.2  & 0.3 \\
        \midrule
        3     & 0.567    & 3.69              & 3.78 & 3.09 & 4.08 & 3.83 \\
        5     & 0.579    & 3.69              & 3.89 & 3.53 & 3.87 & 3.48 \\
        8     & 0.604    & 3.81              & 3.39 & 3.52 & 4.02 & 4.32 \\
        12    & 0.630    & 3.45              & 2.69 & 3.29 & 3.80 & 4.02 \\
        16    & 0.650    & 2.86              & 3.22 & 2.60 & 2.75 & 3.06 \\
        24    & 0.659    & 2.96              & 3.74 & 2.88 & 2.60 & 2.62 \\
        32    & 0.665    & 3.22              & 2.82 & 2.74 & 3.69 & 3.63 \\
        \bottomrule
    \end{tabular}
\end{table}

Within E2b, Table~\ref{tab:e2b_latent_sensitivity} confirms that \(d=16\) provides the lowest aggregate MARE while retaining most of the reconstruction improvement observed over the evaluated dimensions. This value is used for the main E2a, E2b, and E3 experiments. The sensitivity results do not establish that the same dimension is optimal for other geometry families or training regimes.

\section{Input-substitution diagnostics for E3}
\label{app:e3_computational_procedures}

To examine how geometric representations affect the E3 predictions, we use reference-based input substitution following the general baseline-substitution principle used in Shapley-based model interpretation~\cite{chen2022explaining}. The representation input is varied while keeping the trained predictor and Reynolds number fixed.
Each perturbed geometry at \(d_{\mathrm H}>0\) is matched to a reference geometry at \(d_{\mathrm H}=0\) with the same aspect ratio, incident angle, and \(d_{\mathrm L}\).
The matched geometries are independent realizations under these parameter conditions.
Both the original and substituted inputs are evaluated against the same CFD target for the perturbed geometry.
The substitutions therefore test how the information supplied to each trained predictor affects its prediction error, without retraining either model.

The results below focus on the strongest perturbation, \(d_{\mathrm H}=0.3\), comprising 24 matched geometry pairs evaluated at nine Reynolds numbers.
MARE is calculated over the resulting 216 flow observations separately for each training seed.
Both models use training seeds 42--46 with the data split fixed at seed 42; the tables report the mean and sample standard deviation of MARE across these five runs.

\paragraph{Descriptor restoration}
For Descriptor-MLP, we replace one descriptor of the perturbed geometry with its value for the matched reference geometry.
The remaining descriptors retain their perturbed-geometry values.
This gives four input conditions: the original descriptor vector and three variants restoring \(D_{\mathrm{eq}}\), \(\Phi\), or \(\Phi_{\mathrm{cross}}\), respectively.

\begin{table}[htbp]
    \centering
    \caption{Descriptor-MLP input-restoration results at \(d_{\mathrm H}=0.3\). Values are the mean and sample standard deviation across five training seeds. All input conditions are evaluated against the same perturbed-geometry CFD targets.}
    \label{tab:e3_descriptor_restoration}
    \begin{tabular}{lc}
        \toprule
        Input condition                   & MARE (\%) \\
        \midrule
        Original                          & \(9.88\pm3.22\) \\
        Restore \(D_{\mathrm{eq}}\)       & \(9.56\pm3.17\) \\
        Restore \(\Phi\)                  & \(3.28\pm0.50\) \\
        Restore \(\Phi_{\mathrm{cross}}\) & \(4.50\pm1.62\) \\
        \bottomrule
    \end{tabular}
\end{table}

Table~\ref{tab:e3_descriptor_restoration} shows that restoring either area-related descriptor substantially reduces prediction error, whereas restoring \(D_{\mathrm{eq}}\) has little effect.
This pattern holds in all five training runs.
Together with the amplified drag response at the strongest perturbation reported in Section~\ref{sec:e3_extrapolative_analysis}, these results associate the elevated Descriptor-MLP error with its response to changes in \(\Phi\) and \(\Phi_{\mathrm{cross}}\).

\paragraph{Latent-input substitution}
For VAE-MLP, let \(\bm\mu_0\) and \(\bm\mu_1\) denote the encoder posterior means of the matched reference and perturbed geometries, respectively.
We compare the prediction obtained using the complete geometry-specific input \(\bm\mu_1\) with that obtained after substituting the reference input \(\bm\mu_0\).
The trained latent-to-drag predictor and Reynolds number remain unchanged, and both predictions are evaluated against the perturbed-geometry CFD target.

\begin{table}[htbp]
    \centering
    \caption{VAE-MLP latent-input substitution results at \(d_{\mathrm H}=0.3\). Values are the mean and sample standard deviation across five training seeds. Both inputs are evaluated against the same perturbed-geometry CFD targets.}
    \label{tab:e3_latent_substitution}
    \begin{tabular}{lc}
        \toprule
        Input condition                  & MARE (\%) \\
        \midrule
        Matched reference, \(\bm\mu_0\)  & \(6.34\pm0.74\) \\
        Perturbed geometry, \(\bm\mu_1\) & \(4.97\pm0.78\) \\
        \bottomrule
    \end{tabular}
\end{table}

Using the perturbed geometry's latent representation reduces the cross-seed mean MARE from \(6.34\%\) to \(4.97\%\), as shown in Table~\ref{tab:e3_latent_substitution}.
The improvement occurs in every training run, with reductions ranging from \(1.01\) to \(1.73\) percentage points.
The latent representation of the perturbed geometry therefore supplies useful predictive information beyond that provided by the matched reference representation.
This result complements the attenuated cross-level response observed in the main text: the smaller predicted drag changes still incorporate geometry-specific information that improves prediction.


\bibliographystyle{elsarticle-num-names}
\bibliography{references}

\end{document}